# Deterministic Numerical Schemes for the Boltzmann Equation

BY AKIL NARAYAN AND ANDREAS KLÖCKNER <{anaray,kloeckner}@dam.brown.edu>

**Abstract**

This article describes methods for the deterministic simulation of the collisional Boltzmann equation. It presumes that the transport and collision parts of the equation are to be simulated separately in the time domain. Time stepping schemes to achieve the splitting as well as numerical methods for each part of the operator are reviewed, with an emphasis on clearly exposing the challenges posed by the equation as well as their resolution by various schemes.

## 1 Introduction

The Boltzmann equation is an equation of statistical mechanics describing the evolution of a rarefied gas.

- In a *fluid* in *continuum mechanics*, all particles in a spatial volume element are approximated as having the same velocity.

- In a *rarefied gas* in *statistical mechanics*, there is enough space that particles in one spatial volume element may have different velocities.

The equation itself is a nonlinear integro-differential equation which decribes the evolution of the density of particles in a monatomic rarefied gas. Let the density function be $f(x,v,t)$. The quantity $f(x,v,t)\mathrm{d}x\,\mathrm{d}v$ represents the number of particles in the phase-space volume element $\mathrm{d}x\,\mathrm{d}v$ at time $t$. Both $x$ and $v$ are three-dimensional independent variables, so the density function is a map $f\colon \mathbb{R}^3 \times \mathbb{R}^3 \times \mathbb{R} \to \mathbb{R}_0^+ := \mathbb{R} \cup \{0, \infty\}$. The theory of kinetics and statistical mechanics gives rise to the laws which the density function $f$ must obey in the absence of external forces. From these physical constraints, an evolution law for the density can be derived, which reads

$$\frac{\partial f}{\partial t} + v \cdot \nabla_x f = \frac{1}{k_n} Q(f,f), \qquad x, v \in \mathbb{R}^3.$$

The right-hand side terms $k_n$ and $Q(\,\cdot\,,\,\cdot\,)$ are the *Knudsen number* and the *collision operator* and we will discuss them in detail later in this introduction. In a more general model, we should also allow for the existence of external forces: we introduce the $\mathbb{R}^3$-valued function $L(x,t)$, which we can now introduce to form the full *Boltzmann equation*:

$$\frac{\partial f}{\partial t} + v \cdot \nabla_x f + L \cdot \nabla_v f = \frac{1}{k_n} Q(f,f), \qquad x, v \in \mathbb{R}^3. \tag{1}$$

$L$ is often the *Lorentz force* exerted upon charged particles by electromagnetic fields.

The left hand side of (1) is called the *Vlasov Equation*–it is simply advection: the density (i.e. the particles which the density represents) evolves according to its velocity and the Eulerian forces acting upon it. The right hand side operator $Q(f,f)$ is a *collision term*, which represents the binary collisions experienced by individual particles. The scalar $k_n$ (the Knudsen number) is a dimensionless quantity that is

defined as follows

$$k_n := \frac{\lambda}{l} := \frac{\text{mean free path}}{\text{representative physical length scale}}.$$

In all the following, we shall adopt the notation presented in [15], with small changes. The collision operator $Q(\,\cdot\,,\,\cdot\,)$ is derived from physical considerations and can be broken up into a gain and a loss term:

$$Q(f,f) = Q^+(f,f) - L[f]f, \tag{2}$$

where the gain and loss terms are defined by

$$Q^+(f,f) = \int_{\mathbb{R}^3} \int_{S^2} B(|v-v_*|, \cos\theta)\, f(v')\, f(v'_*)\, \mathrm{d}\omega\, \mathrm{d}v_*, \tag{3}$$

$$L[f] = \int_{\mathbb{R}^3} \int_{S^2} B(|v-v_*|, \cos\theta)\, f(v_*)\, \mathrm{d}\omega\, \mathrm{d}v_*. \tag{4}$$

$B(\,\cdot\,,\,\cdot\,)$ is a collision kernel modeling the mechanics of a collision, and $S^2 \subset \mathbb{R}^3$ is the unit sphere. Figure 1 illustrates the meaning of these variables. In particular, observe that $\omega$ is the (unit-length) pre-collision deflection relative to path of the center of gravity of the system consisting of the two particles.

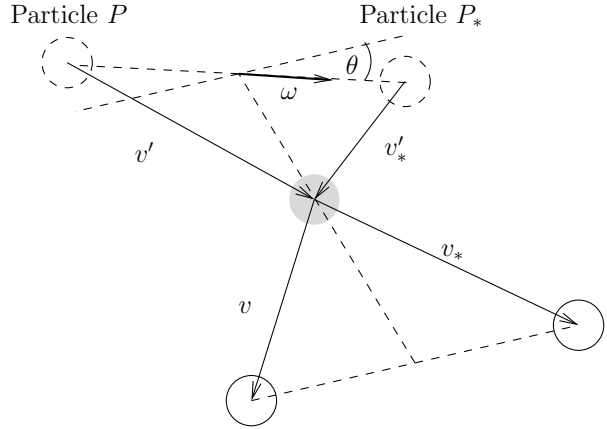

**Figure 1.** Illustration of the notation in (3) and (4).

The quantities $v$ and $v_*$ represent the velocities of two particles post-collision, and $v'$ and $v'_*$ are the same quantities pre-collision. Given the post-collision velocities and a deflection vector $\omega$, the pre-collision velocities may be calculated by

$$v' = \frac{1}{2}(v + v_* + |v - v_*|\omega), \qquad v'_* = \frac{1}{2}(v + v_* - |v - v_*|\omega).$$

The angle $\theta$ is the angle between $v - v_*$ and $v' - v'_*$. To gain an intuitive insight, consider the special cases $\theta = 0°$ (head-on collision) and $\theta = 90°$ (one particle in motion, the other at rest). In particular,

$$\cos\theta = \frac{\omega \cdot (v - v_*)}{|v - v_*|}.$$

On the gain side, the term

$$B(|v-v_*|, \cos\theta)\, \underbrace{f(v')\, f(v'_*)}_{\text{availability}}\mathrm{d}\omega\, \mathrm{d}v_*$$

describes, for given post-collision velocities $v$ and $v_*$ and given deflection vector $\omega$, "how many" particles (relatively) take part in this collision and acquire the new velocity $v$. Of course, the gain depends on "how many" particles of the pre-collision velocities are available, hence we multiply by the pre-collision densities

$$f(v')f(v'_*)$$

highlighted as "availability" above.

Similarly the loss term

$$B(|v-v_*|,\cos\theta)\underbrace{f(v)f(v_*)}_{\text{availability}}\mathrm{d}\omega\,\mathrm{d}v_*$$

describes, this time for given *pre*-collision velocities $v$ and $v_*$ and given deflection vector $\omega$, how many particles leave the current velocity state. Of couse, we again need to multiply by an availability factor, as highlighted above. Note that each of these considerations needs to be repeated at *each* point $x$ in space. Several examples (in order of increasing complexity) of kernels $B$ are

- *Maxwellian gas:* $B(|v-v_*|,\cos\theta) = \text{const}$.
- *Hard Sphere gas:* $B(|v-v_*|,\cos\theta) = \text{const}\cdot|v-v_*|$.
- *Variable Hard Sphere (VHS) gas:* $B(|v-v_*|,\cos\theta) = \text{const}\cdot|v-v_*|^\alpha$. (generalizes both cases above)

Solving the Boltzmann equation is tantamount to:

- solving the 6+1, space+time, nonlinear integrodifferential equation given by (1),
- computing the collision operator (2) by performing the integrals in (3) and (4).

In this survey, we shall only discuss methods to deal with the first two items: the Boltzmann equation (1) and the collision term (2). We will assume that the forcing function $L$ can be readily computed by some means.

Even if we assume the force $L$ is available cheaply, one can see that the Boltzmann equation is certainly one of the harder equations in scientific computing. In particular, it presents the following challenges:

- *Dimensionality*: The density function $f$ is defined on a six-dimensional space; three of those dimensions are unbounded. The collision operator contains another five-dimensional integral–for each point in six-dimensional space. A naive program would have to contain eleven nested loops–one over each dimension. For comparison, modern supercomputers already have difficulty with three-dimensional problems at moderate resolutions.
- *Timestepping*: The collision operator and the advection require fairly different approaches to timestepping, as we will discuss later. How can the two different methods be merged?
- *Collision Operator*: While we cannot dodge the six dimensions of the space that $f$ lives on, we need to be smart about carrying out the collision integrals. How can that be done?
- *Shock waves and Positivity for the Transport*: While the transport part of the equation is not the main consumer of computational time, it is still far from straightforward in that it must maintain positivity and resolve shock waves, while efficiently dealing with the high dimensionality of the data.

In this survey, we shall summarize some of the recent methods that have been used to overcome these difficulties and weigh the advantages and disadvantages of each. We will deal excluusively with deterministic methods, noting that these are best applied near an equilibrium state, when additional accuracy is desired. In non-equilibrium situations, or under less stringent accuracy requirements, particle methods may be more feasible.

We organize this survey as follows: We start by describing some methods for time-stepping in Section 2; these methods form the foundation of the splitting strategy, which is at the heart of every method presented here. In Section 3, we describe how one may handle the linear advection components of the Boltzmann equation (i.e. the Vlasov equation). In Section 4, we describe some methods to efficiently compute the collision operator (2). Each section is designed to be meaningfully read on its own.

## 2 Splitting Methods

The full Boltzmann equation (1) is a very challenging equation to solve: the transport and collision contributions must be solved in very different ways (see Sections 3 and 4). Because of this, we shall seek ways to solve the transport contribution and the collision contribution individually. Once we have developed robust high-order methods to solve each part, we shall attempt to fuse them together to form a complete scheme for the full Boltzmann equation. The numerical methods which serve as the algorithmic glue of the individual schemes are the time-stepping procedures described in this section.

In Sections 3 and 4 we will describe methods which are used to solve two separate PDE systems. In Section 3 describing the transport equation, we shall be concerned with solving

$$\frac{\partial f}{\partial t} = -v \cdot \nabla_x f - L \cdot \nabla_v f =: \mathcal{A}[f]. \tag{5}$$

Similarly, in Section 4 pertaining to collision we want to solve

$$\frac{\partial f}{\partial t} = \frac{1}{k_n} Q(f, f) =: \mathcal{B}[f]. \tag{6}$$

However, we do not wish to solve (5) and (6) separately. Instead, we wish to solve the full Boltzmann equation (1), which we can rewrite as

$$\frac{\partial f}{\partial t} = \mathcal{A}[f] + \mathcal{B}[f]. \tag{7}$$

Unfortunately there is no simple way to solve (7) by solving (5) and (6) alone: suppose $f_1$ solves (5) and $f_2$ solves (6). Then unless both $\mathcal{A}$ and $\mathcal{B}$ are linear and $\mathcal{A}[f_2] = -\mathcal{B}[f_1]$, then $f_1 + f_2$ does not solve (7). For the Boltzmann equation, $\mathcal{A}$ is not linear due to the Lorentz force $L$, and $Q$ is also nonlinear. Even if we ignore the nonlinearity, we cannot expect to satisfy $\mathcal{A}[f_2] = -\mathcal{B}[f_1]$. Thus superposition is not a viable solution method.

Suppose instead that we just derive semi-discrete formulations for (5) and (6) separately and then combine the resulting right-hand-side phase-space discretizations. To be precise, consider a transport method, e.g. one of the methods discussed in Section 3 which has degrees of freedom $\{f_i^{\mathcal{A}}\}_{i \in I}$ and the semi-discrete form

$$\frac{\mathrm{d} f_i^{\mathcal{A}}}{\mathrm{d} t} = \mathcal{A}_h[f_i^{\mathcal{A}}] \quad \forall i \in I, \tag{8}$$

for some discretization of the transport operator $\mathcal{A}_h$. Consider also a collision method chosen from e.g. Section 4 with degrees of freedom $\{f_j^{\mathcal{B}}\}_{j \in J}$ accompanied by the assoicated semi-discrete form

$$\frac{\mathrm{d} f_j^{\mathcal{B}}}{\mathrm{d} t} = \mathcal{B}_h[f_j^{\mathcal{B}}] \quad \forall j \in J. \tag{9}$$

One might suggest somehow melding these methods into one for the degrees of freedom $\{f_k\}_{k\in K}$ for the full Boltzmann equation (7) which reads

$$\frac{\mathrm{d}f_k}{\mathrm{d}t} = \mathcal{P}_i^k \mathcal{A}_h[f_i^\mathcal{A}] + \mathcal{P}_j^k \mathcal{B}_h[f_j^\mathcal{B}] \quad \forall k \in K, \tag{10}$$

where $\mathcal{P}_i^k$ is some transformation operator which maps the degrees of freedom $f_i^\mathcal{A}$ into the degrees of freedom $f_k$, and $\mathcal{P}_j^k$ is defined analogously. Indeed this is quite tempting until one realizes that there exists a great many problems with this formulation. Firstly, what are the transformation operators $\mathcal{P}_i^k$, $\mathcal{P}_j^k$ and their inverses $\mathcal{P}_k^i$, $\mathcal{P}_k^j$? If the methods do not share the same type of degrees of freedom (e.g. the transport discretization is discontinuous Galerkin and the collision scheme is a Fourier method), then transferring between the $f_k$, and the $f_i^\mathcal{A}$ and $f_j^\mathcal{B}$ may prove complicated and/or onerous.

In addition, one quickly realizes that many of the numerical methods presented in Section 3 do not have a semi-discrete formulation (8). Indeed, many of the methods presented in that section integrate the degrees of freedom exactly in time: there is no semi-discrete form. Methods like this cannot be molded into forms (8). Thus, the formulation (10) fails for many of the methods described in Section 3.

Another drawback is that even if the transport and the collision degrees of freedom are the same (i.e. $f_i^\mathcal{A}$, $f_j^\mathcal{B}$, and $f_k$ are all the same degrees of freedom), using one single explicit time-stepping procedure on (10) must adhere to the strictest stability criterion in both $\mathcal{A}_h$ and $\mathcal{B}_h$. We might then suggest using an implicit method for (9) to avoid the timestep restriction, while using a high-order explicit method on (8) to maintain good accuracy for the transport term. However, using two separate time-stepping methods is not possible in a straightforward implementation of (10).

This leads us to the topic of this section: time-splitting methods. We suppose that the degrees of freedom $f_i^\mathcal{A}$ and $f_j^\mathcal{B}$ are the same type of discretization so that we may simply write $f_k^\mathcal{A}$ and $f_k^\mathcal{B}$ as the degrees of freedom. We also introduce a superscript indexed by $n$, the time step; i.e. $f_k^n$ is an approximation to $f_k(t^n)$, where the time $t$ is indexed by $n$. We now suppose that we wish to solve

$$\frac{\mathrm{d}f_k}{\mathrm{d}t} = \mathcal{A}[f_k] + \mathcal{B}[f_k], \tag{11}$$

where we have dropped the subscript $h$ on the operators $\mathcal{A}$ and $\mathcal{B}$ for readability. The augmented variables $f_k^\mathcal{A}$ and $f_k^\mathcal{B}$ are defined according to the evolution laws

$$\left.\begin{aligned}\frac{\mathrm{d}f_k^\mathcal{A}}{\mathrm{d}t} &= \mathcal{A}[f_k^\mathcal{A}], \\ \frac{\mathrm{d}f_k^\mathcal{B}}{\mathrm{d}t} &= \mathcal{B}[f_k^\mathcal{B}].\end{aligned}\right\} \tag{12}$$

Suppose we choose a separate time discretization for each part of (12), and let $S_\mathcal{A}$ and $S_\mathcal{B}$ be the solution operators for the differential equations in (12). We can rewrite (12) as

$$f_k^\mathcal{A}(t^{n+1}) = S_\mathcal{A}(f_k^\mathcal{A}(t^n); t^n, t^{n+1}),$$

$$f_k^\mathcal{B}(t^{n+1}) = S_\mathcal{B}(f_k^\mathcal{B}(t^n); t^n, t^{n+1}).$$

If the semi-discrete form is time-autonomous, then we simply take $S$ as an operator with a superscript denoting the time step $t^{n+1} - t^n =: \Delta t$ and write

$$\left.\begin{aligned}f_k^\mathcal{A}(t^{n+1}) &= S_\mathcal{A}^{\Delta t} f_k^\mathcal{A}(t^n), \\ f_k^\mathcal{B}(t^{n+1}) &= S_\mathcal{B}^{\Delta t} f_k^\mathcal{B}(t^n).\end{aligned}\right\} \tag{13}$$

In addition, we introduce stepping operators $S_{\mathcal{A},n}^{\Delta t}$ and $S_{\mathcal{B},n}^{\Delta t}$ as the temporal discretizations associated with (12). Then we can write

$$\left.\begin{aligned} f_k^{\mathcal{A},n+1} &= S_{\mathcal{A},n}^{\Delta t} f_k^{\mathcal{A},n}, \\ f_k^{\mathcal{B},n+1} &= S_{\mathcal{B},n}^{\Delta t} f_k^{\mathcal{B},n}. \end{aligned}\right\} \tag{14}$$

The stepping operators do not have to be associated to any temporal discretization. For example, $S_{\mathcal{A},n}^{\Delta t}$ may be given by the flux balance method advancement in equation (35). The advantage of this formulation is that we have written the system (5) and (6) in both semi-discrete form (13) and fully discrete form (14). We shall use both the semi-discrete and the fully discrete formulations to describe our time-stepping procedures.

In the following sections, we shall assume that we are trying to solve the full system (7) so that the superscript $\mathcal{A}$ or $\mathcal{B}$ on the unknowns $f_k$ becomes unnecessary as we are no longer solving two different equations (5) and (6). Thus, we shall write the semidiscrete degrees of freedom as $f_k(t)$ and the fully discrete ones as $f_k^n$.

## 2.1 Operator Splitting

Our first solution to the splitting problem is the operator splitting assertion. Assume we have a linear ODE

$$\frac{\mathrm{d}y}{\mathrm{d}t} = (A+B)y,$$

where $A$ and $B$ are some operators. The solution to this system can be written as

$$y = e^{(A+B)t}.$$

We wish to deal with $A$ and $B$ separately, so we may attempt to form a scheme which tries to approximate

$$e^{(A+B)t} \approx P(e^{At}, e^{Bt}),$$

for some function $P$. That is, instead of computing $e^{(A+B)t}$, we instead compute $e^{At}$ and $e^{Bt}$ and combine them in some clever fashion so that we still obtain temporal accuracy and convergence. We now survey some of the most popular methods of this type that are used.

### 2.1.1 Simple First-Order

The simplest thing one could think of to employ both of the stepping operators $S_{\mathcal{A}}$ and $S_{\mathcal{B}}$ is to use them in a sequential manner. I.e.

$$f_k(t^{n+1}) = S_{\mathcal{B}}^{\Delta t} S_{\mathcal{A}}^{\Delta t} f_k(t^n). \tag{15}$$

This amounts to solving the following ODE systems in a sequential manner at each time level:

$$\left.\begin{aligned} \frac{\mathrm{d}\tilde{f}_k(t)}{\mathrm{d}t} &= \mathcal{A}[\tilde{f}_k(t)], \\ \tilde{f}_k(t^n) &= f_k^n, \end{aligned}\right\} \tag{16}$$

$$\left.\begin{aligned} \frac{\mathrm{d}f_k(t)}{\mathrm{d}t} &= \mathcal{B}[f_k(t)], \\ f_k(t^n) &= \tilde{f}_k(t^{n+1}). \end{aligned}\right\} \tag{17}$$

We can then set calculate $f_k(t^{n+1})$ as given by the solution to system (17). The advantage here is clear: we can independently apply methods from sections 3 and 4 to solve (16) and (17), respectively. Thus, all the analysis and validation from those sections apply here. The downside is that one can show that this simple method is only first-order accurate in time (see e.g. [38]).

### 2.1.2 Second-Order Strang Splitting

A very popular alternative to the suboptimal first-order splitting technique presented in the previous section is a splitting strategy first suggested by Strang in [48]. The main idea is to symmetrically straddle the operators on the temporal intervals. I.e.

$$f_k(t^{n+1}) = S_{\mathcal{A}}^{\Delta t/2} S_{\mathcal{B}}^{\Delta t} S_{\mathcal{A}}^{\Delta t/2} f_k(t^n), \tag{18}$$

Strang shows that this method is second-order accurate. Thus, as long as the fully discrete operators $S_{\mathcal{A},n}$ and $S_{\mathcal{B},h}$ are second-order accurate then the entire scheme will be second-order accurate in time. The Strang splitting technique is one of the most popular algorithms for solving this type of problem due to its simplicity and its relative accuracy. Figure 2 visually compares the methodologies of first-order and Strang operator splitting.

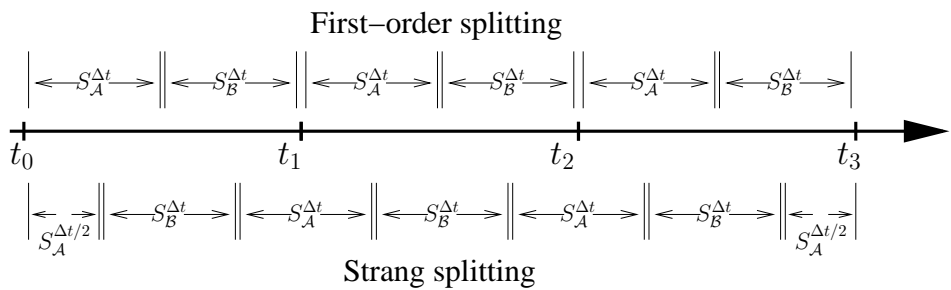

**Figure 2.** An illustration of the first-order splitting method and the Strang splittng method.

### 2.1.3 More complex strategies

The Strang and first-order splitting methods are the most commonly employed methods. There are other methods available, but these methods are more complicated and they do not provide any extraordinary accuracy. For example, there is another second-order splitting method in constrast to Strang's method due to Ohwada [38]. For this method, we explicity define the fully-discrete operator

$$S_{\mathcal{B},n}^{\Delta t} f_k^n = f_k^n + \Delta t\, \mathcal{B}[f_k^n],$$

which is basically a forward-Euler method. This allows us to define the scheme as

$$\begin{aligned} f_k^{n+1} &= S_{\mathcal{A}}^{\Delta t} f_k^n + \frac{1}{2}\bigl(S_{\mathcal{B},n}^{\Delta t} S_{\mathcal{B},n}^{\Delta t} S_{\mathcal{A}}^{\Delta t} f_k^n - S_{\mathcal{B},n}^{\Delta t} S_{\mathcal{A}}^{\Delta t} f_k^n\bigr) + \frac{1}{2}\bigl(S_{\mathcal{B},n}^{\Delta t} S_{\mathcal{A}}^{\Delta t} f_k^n - S_{\mathcal{A}}^{\Delta t} f_k^n\bigr) \\ &= \frac{1}{2}\bigl(\mathcal{I} + S_{\mathcal{B},n}^{\Delta t} S_{\mathcal{B},n}^{\Delta t}\bigr) S_{\mathcal{A}}^{\Delta t} f_k^n. \end{aligned} \tag{19}$$

Ohwada shows that this scheme is second-order accurate, and that it preserves mass, momentum, and energy in the context of the Boltzmann equation. (See Section 4.1 for mathematical definitions of these physical quantities.) Note that this fully-discrete scheme uses the exact solution operator $S_{\mathcal{A}}^{\Delta t}$ and not some temporal discretization of it. However, given that $\mathcal{A}$ is simply convection, this can be solved exactly in principle. The Ohwada splitting scheme (19) has roughly the same computational complexity as the Strang splitting, with the advantage of conservation of relevant quantities for the Boltzmann equation.

There are also other types of higher-order splitting methods in e.g. [13]. For example, we can can consider two different operator splitting schemes of the form

$$f_k(t^{n+1}) = \frac{1}{3}\Big( 4\, S_{\mathcal{A}}^{\Delta t/4} S_{\mathcal{B}}^{\Delta t/2} S_{\mathcal{A}}^{\Delta t/2} S_{\mathcal{B}}^{\Delta t/2} S_{\mathcal{A}}^{\Delta t/4} - S_{\mathcal{A}}^{\Delta t/2} S_{\mathcal{B}}^{\Delta t} S_{\mathcal{A}}^{\Delta t/2} \Big) f_k(t^n), \tag{20}$$

$$f_k(t^{n+1}) = \frac{1}{6}\Big( 4\, S_{\mathcal{A}}^{\Delta t/2} S_{\mathcal{B}}^{\Delta t} S_{\mathcal{A}}^{\Delta t/2} + 4 S_{\mathcal{B}}^{\Delta t/2} S_{\mathcal{A}}^{\Delta t} S_{\mathcal{B}}^{\Delta t/2} - 2\, S_{\mathcal{A}}^{\Delta t} S_{\mathcal{B}}^{\Delta t} - S_{\mathcal{B}}^{\Delta t} S_{\mathcal{A}}^{\Delta t} \Big) f_k(t^n). \tag{21}$$

Dia [13] shows that these splitting strategies are fourth-, and third-order accurate in time, respectively. However, they suffer from some stability issues and thus are not used very extensively. In addition, evaluating $S_{\mathcal{B}}^{\Delta t}$ (the collision operator evolution) multiple times can be a computational burden.

## 2.2 Implicit-Explicit Methods

To avoid using the operator splitting strategy, implicit-explicit Runge-Kutta (IMEX-RK) schemes have been devised which appear to boast the advantages of the splitting strategy without the unfortunate restriction to second-order. These methods were first introduced by Ascher et al. in [2]. The basic idea of these methods is as follows: recall that we have an ODE system

$$\frac{\mathrm{d} f_k}{\mathrm{d} t} = \mathcal{A}_h[f_k] + \mathcal{B}_h[f_k]. \tag{22}$$

The operator-splitting approach sought to find solution operators for each right-hand side contribution individually and then to apply them in a clever way to attain consistency and accuracy. The IMEX-RK approach is exactly the same, except that we shall introduce a very general mathematical framework by which we can explicitly derive schemes. In order to conform to the notation in the literature, we rewrite (22) in the form

$$y' = f(y) + g(y), \tag{23}$$

where $f(y)$ is some non-stiff, possibly highly nonlinear, term and $g(y)$ is a stiff, but 'close' to linear, term. The stiffness of $g(y)$ can be attributed to implicit dependence on some relaxation parameter $\varepsilon$; we shall not discuss the effect of the degree of stiffness on the numerical solution of (23), but in general, we should replace $g(y)$ with $\frac{1}{\varepsilon} g(y)$, and consider stability and accuracy as functions of $\varepsilon$, which controls the stiffness of the system. However, we shall make only brief mentions of the relaxation parameter $\varepsilon$ in the future, and so omit its explicit role in equation (23). For the Boltzmann equation, the relaxation parameter $\varepsilon$ is the Knudsen number, previously labelled $k_n$.

For example, (23) may arise from some convection-diffusion PDE, with $f(y)$ representing the semi-discrete form of the convection, and $g(y)$ representing the semi-discrete form of the diffusion. In this way, diffusive second-derivative terms are more stiff than convective terms. In the context of the Boltzmann equation, $f(y)$ represents the convective Vlasov equation terms, and $g(y)$ represents the possibly stiff collision operator. Note that we have, without loss, suppressed any non-autonomous dependence in (23): a non-autonomous system can always be written as an autonomous one, but in any case we shall develop the IMEX-RK methods with non-autonomous generality. In the following, we write $f = f(t, y)$ and $g = g(t, y)$.

An IMEX-RK scheme is defined by the number of *stages* $s$; by four $s \times 1$ vectors $b$, $\tilde{b}$, $c$, and $\tilde{c}$; and by two $s \times s$ matrices $a$ and $\tilde{a}$. The method takes the unknowns $y^n$ at time level $t^n$ to time level $t^{n+1} = t^n +$

$\Delta t$ via the scheme

$$y^{n+1} = y^n + \Delta t \sum_{j=1}^{s} \tilde{b}_j f(t^n + \tilde{c}_j \Delta t, Y_j) + \Delta t \sum_{j=1}^{s} b_j g(t^n + c_j \Delta t, Y_j), \tag{24}$$

where the intermediate stage variables $\{Y_i\}_{i=1}^{s}$ are given by

$$Y_i = y^n + \Delta t \sum_{j=1}^{s} \tilde{a}_{ij} f(t^n + \tilde{c}_j \Delta t, Y_j) + \Delta t \sum_{j=1}^{s} a_{ij} g(t^n + c_j \Delta t, Y_j). \tag{25}$$

(24) and (25) form the IMEX-RK scheme. So far, we have only stated a generalization of the usual RK scheme. Indeed, if $b = \tilde{b}$, $c = \tilde{c}$, and $a = \tilde{a}$, then the IMEX-RK scheme degenerates into a familiar $s$-stage RK scheme. The weights $a$, $b$, and $c$, and their tilde'd counterparts are usually arranged in a *Butcher tableau* [9] as depicted in Figure 3.

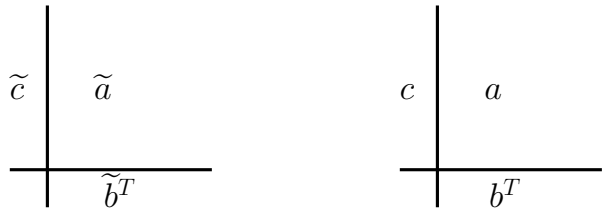

**Figure 3.** Arrangement of Butcher tableaux.

Recall that our motivation for this particular type of splitting is that the term $f$ is usually relatively straightforward to integrate with an explicit method. Thus, we shall define the tableau for $f$ (the tilde'd variables) to be an explicit RK method. Recall that in the arrangement of the tableau, the nonzero entries of the matrix $\tilde{a}$ determine the implicit or explicit nature of a method:

- An *explicit* method has zero entries on and above the main diagonal: $\tilde{a}_{ij} = 0$ for $j \geq i$.

- A *diagonally implicit* method has zeros above the main diagonal: $\tilde{a}_{ij} = 0$ for $j > i$.

- A method not satisfying any of the above two restrictions is called fully *implicit*.

We shall restrict our attention to solving $f$ via an explicit method. Thus, we hereafter assume that $\tilde{a}_{ij} = 0$ for $j \geq i$. In addition, although we can allow $g$ to be solved with any fully implicit method, practical considerations may become problematic if $g$ is not strictly linear. Thus, we also assume that $g$ is solved with a diagonally implicit Runge-Kutta (DIRK) method, which satisfies $a_{ij} = 0$ for $j > i$.

The purpose of this survey is not to describe the theory of IMEX-RK methods. All that we shall say is that in order to determine the order of accuracy for traditional RK schemes, the essential method is Taylor expansion: the exact solution is substituted into the numerical scheme, is Taylor expanded, and the highest-magnitude Taylor expansion term which is not eliminated by the scheme dictates the order of convergence in $\Delta t$. However, these IMEX-RK methods are not traditional RK methods; in fact, one can view these split methods as special cases of *partitioned* RK methods [28]. Thus, we can apply the order condition derivation presented by Hairer in [25] to determine the order of convergence of these IMEX-RK methods. In fact, Pareschi and Russo derive some of these conditions in [40], which we present below to illustrate the characteristics of the method.

### 2.2.1 Order of Convergence

Before providing some examples of these methods, we shall show and interpret some of the order conditions for IMEX-RK schemes that have been derived in [40]. We suppose that we are looking for a scheme of order $p$ in $\Delta t$. This requirement leads to a set of conditions on the values of the elements in the Butcher tableau. We present conditions for the first three orders below:

- $p=1$. The elements in the Butcher tableau in Figure 3 must satisfy

$$\sum_{i=1}^{s} \tilde{b}_i = 1, \quad \sum_{i=1}^{s} b_i = 1.$$

- $p=2$. The elements in the Butcher tableau must satisfy

$$\sum_{i=1}^{s} \tilde{b}_i \tilde{c}_i = \frac{1}{2}, \quad \sum_{i=1}^{s} b_i c_i = \frac{1}{2},$$

$$\sum_{i=1}^{s} \tilde{b}_i c_i = \frac{1}{2}, \quad \sum_{i=1}^{s} b_i \tilde{c}_i = \frac{1}{2}.$$

- $p=3$. The elements in the Butcher tableau must satisfy

$$\sum_{i,j=1}^{s} \tilde{b}_i \tilde{a}_{ij} \tilde{c}_j = \frac{1}{6}, \quad \sum_{i,j=1}^{s} b_i a_{ij} c_j = \frac{1}{6}, \quad \sum_{i=1}^{s} \tilde{b}_i \tilde{c}_i \tilde{c}_i = \frac{1}{3}, \quad \sum_{i=1}^{s} b_i c_i c_i = \frac{1}{3},$$

$$\sum_{i,j=1}^{s} \tilde{b}_i a_{ij} c_j = \frac{1}{6}, \quad \sum_{i,j=1}^{s} b_i \tilde{a}_{ij} c_j = \frac{1}{6}, \quad \sum_{i,j=1}^{s} b_i a_{ij} \tilde{c}_j = \frac{1}{6}, \quad \sum_{i,j=1}^{s} \tilde{b}_i \tilde{a}_{ij} c_j = \frac{1}{6},$$

$$\sum_{i,j=1}^{s} b_i \tilde{a}_{ij} \tilde{c}_j = \frac{1}{6}, \quad \sum_{i,j=1}^{s} \tilde{b}_i a_{ij} \tilde{c}_j = \frac{1}{6},$$

$$\sum_{i=1}^{s} \tilde{b}_i c_i c_i = \frac{1}{3}, \quad \sum_{i=1}^{s} b_i \tilde{c}_i c_i = \frac{1}{3}, \quad \sum_{i=1}^{s} \tilde{b}_i \tilde{c}_i c_i = \frac{1}{3}, \quad \sum_{i=1}^{s} b_i \tilde{c}_i \tilde{c}_i = \frac{1}{3}.$$

We can observe some very interesting properties from the order conditions presented above: for first-order approximations consistency of the entire method is guaranteed by consistency of the individual methods. I.e., the requirement for first-order consistency places no restriction on the RK interaction between the tilde'd and non-tilde'd schemes. For second order, the first line of equations is simply the usual second-order consistency requirements for RK schemes; the second line of equations are the consistency requirements because of the coupling of the schemes. The same can be seen with third-order $p=3$: the first line of equations are the natural third-order RK consistency requirements, whereas all the subsequent requirements arise due to the implicit-explicit scheme coupling.

Of course, the error term in deriving these schemes depends heavily on the stiffness in the term $g$ (i.e. the herein omitted relaxation parameter $\varepsilon$). In order to control stiffness of the system, Pareschi and Russo [40] show that $L$-stability of the implicit RK scheme guarantees the advertised order of convergence even for very stiff terms.

### 2.2.2 Examples

In this section we shall present some examples of IMEX-RK methods. Let the triplet $(s, \sigma, p)$ represent the number of implicit stages, number of explicit stages, and the convergence order, respectively. (Note that in our presentation, the number of implicit and explicit stages are the same, but one may introduce intermediate implicit-only stages and then pad the explicit tableau with rows of zeros.) Figure 4 shows two relatively simple IMEX-RK methods: the first is a standard splitting method in which forward-Euler is performed on the explicit part $f$ and backward-Euler is performed on the implicit part $g$. Application of (25) and (24) show that the standard splitting (1,1,1) scheme is equivalent to first-order operator splitting given in (15) with $S_{\mathcal{A},n}$ being forward-Euler and $S_{\mathcal{B},n}$ being backward-Euler.

$$\begin{array}{c|c} 0 & 0 \\ \hline & 1 \end{array} \qquad \begin{array}{c|c} 1 & 1 \\ \hline & 1 \end{array}$$

Standard splitting (1,1,1)

$$\begin{array}{c|cc} 0 & 0 & 0 \\ 1 & 1 & 0 \\ \hline & \frac{1}{2} & \frac{1}{2} \end{array} \qquad \begin{array}{c|cc} -1 & -1 & 0 \\ 2 & 1 & 1 \\ \hline & \frac{1}{2} & \frac{1}{2} \end{array}$$

Jin splitting (2,2,2)

**Figure 4.** Two examples of splitting methods. Explicit tableaus (left) and implicit tableaus (right).

In Figure 4 we also show a second order scheme due to Jin [31]. For many higher-order schemes and further analysis, we refer the reader to [40] and [7].

Many general single-step splitting schemes can be formulated as IMEX-RK schemes. Simple first-order splitting (15) and Strang splitting (18) can be re-cast as IMEX-RK methods. In addition, Jin's (2,2,2) scheme was actually first proposed as an operator-splitting type method.

## 2.3 Other Methods

There have been a great number of other methods which have been proposed; some for solving ODE's of the mixed type (23), some specifically designed for the Boltzmann equation.

Gerisch and Weiner [21] derive some explicit RK methods which they can show preserves positivity of the solution. Of course, this is important in solving the Boltzmann equation which has a non-negative density function as the unknown. The downside is that these methods do not cope well with stiff terms.

If we consider (23) as a function of the relaxation parameter $\varepsilon$, then we can consider what should happen in the limit as $\varepsilon \to \infty$, i.e. as the ODE system becomes a semi-discrete form of a hyperbolic conservation law. It is reasonable to ask that in this limit we obtain some temporal discretization which preserves some desired properties when simulating a hyperbolic conservation law, e.g. the strong stability-preserving (SSP) property. (Formerly called the total variation diminishing property [23].) Such schemes, which are IMEX-RK schemes (and hence retain their robustness in the stiff regime under the assumption of $L$-stability) and preserve the SSP property in the limit of the hyperbolic conservation law, do exist. They are described in a later paper by Pareschi and Russo [42].

We have only described *single-stage* and *multi-stage* methods in this survey, but there are also *multi-step* methods, which require accessing the value of the unknown at previous time-steps. These methods have IMEX variations [3] and SSP-IMEX variations [22].

As a final remark, we have only considered splitting-type algorithms for solving ODEs. The motivation for this stems from the complex nature of the Boltzmann equation. Of course one may use any of the plethora of ODE existing un-split schemes in order to time-step the Boltzmann equation.

We have presented operator splitting methods and IMEX-RK methods, and their variants. Having now described how to connect split transport and collision steps in a time-discretization, we now turn to the problem of actually solving the PDEs associated with each of those operators individually.

# 3 Transport

In this section, we focus on solving a transport equation of the form

$$\frac{\partial f}{\partial t} + v \cdot \nabla_x f + L \cdot \nabla_v f = 0, \tag{26}$$

where $f$, the density function, depends on seven independent variables (including time), and $L$ is the spatially-dependent forcing function (e.g. the electromagnetic Lorentz force). This equation models the natural inertia of the system with external forcing ($L$), and is called the collisionless Boltzmann equation, or more frequently, the *Vlasov equation*. It is an ordinary advection equation, but we must solve it in the full 6-dimensional $(x, v)$ phase space. The equation is implicity nonlinear since the forcing function $L$ usually depend on the density $f$.

In Section 2 we described many methods of time-stepping once a semi-discrete or fully discrete form of (26) has been found. This survey cannot do deserved justice to the voluminous literature on the subject of finding a semi-discrete form for a nonlinear hyperbolic conservation law. ($v$ is $x$-divergence free, and $L$ is $v$-divergence free, so (26) is easily recast into a conservation law form with six 'spatial' dimensions, $x$-$v$ phase space.) We can only mention a brief selection of these methods and point to a few references:

Conservation laws may be converted into semi-discrete form with finite-difference [24], finite-volume [34], finite-element [8], or spectral [26] degrees of freedom. Discontinuous Galerkin methods [12, 27] and nonlinear reconstruction techniques like ENO/WENO [46] provide robust, high-order solvers for nonlinear problems.

Once a semi-discrete method has been found, one can use this in tandem with a solver of choice for the collision term (see Section 4) and apply e.g. an IMEX-RK method, as described in Section 2. One may also choose to adopt an operator-splitting technique, but this requires some choice of time discretization. For methods which have an appropriate semi-discrete form, the transport step has a fairly well-defined shape.

Thus, instead of focusing on well-explored ways to turn (26) into a semi-discrete scheme non-specific to the Boltzmann equation, we shall concentrate on:

- fully-discrete schemes,
- semi-discrete schemes specifically tailored for Vlasov-type equations.

Among the fully-discrete methods we shall consider are the semi-Lagrangian method (discussed in [4] and [17]) and the flux balance method [14]. Some semi-discrete type schemes tailored for the Vlasov equation include Gamba and Proft's DG scheme [20] and Klimas's filtered spectral method [32].

## 3.1 Desired Scheme Properties

The density function is a physical quantity, and we therefore ask our numerical approximation of the density to satisfy certain properties:

- Positivity: $f(x, v, t) \geq 0 \ \forall \ x, v, t$
- Conservation
  - Mass: $\int f \, dv \, dx = m$ for some constant scalar $m$ independent of time
  - Momentum: $\int v f \, dv \, dx = M$ for some constant vector $M$ independent of time
  - Energy: $\int |v|^2 f \, dv \, dx = E$ for some constant scalar $E$ independent of time

Unfortunately, these restrictions (especially positivity) prevent us from blindly applying the usual well-developed numerical schemes for hyperbolic conservation laws. In addition, conservation of momentum and energy is an extremely difficult task and few numerical methods we present here actually address these latter two goals. The following sections describe some of the high-order accurate numerical methods that have been presented to deal with these difficulties.

## 3.2 More Splitting

Solving (26) requires a six-dimensional solver which can be quite complex to formulate and code. Because of this, many have observed that we can take advantage of operator splitting to simplify our approach. We can split the physical space $x$-advection and the velocity space $v$-advection. We obtain the collection of equations

$$\left.\begin{aligned} \frac{\partial f}{\partial t} + v \cdot \nabla_x f &= 0, \\ \frac{\partial f}{\partial t} + L \cdot \nabla_v f &= 0. \end{aligned}\right\} \qquad (27)$$

Since $v$ is independent of $x$ and $L$ is independent of $v$, then each of these equations is a scalar multidimensional hyperbolic equation of the form

$$\frac{\partial u}{\partial t} + c \cdot \nabla u = 0. \qquad (28)$$

That is, we are dealing with constant-coefficient advection. (Note that $L$ is not really independent of $f$, so it is not strictly linear; however, the dependence on $f$ tends to be complicated, and so we shall consider $L$ a given function.) We can then solve each equation in (27) separately and then use an operator splitting approach to advance the density in time.

Note that even devising methods for solving (28) is not trivial since we ask for positivity and conservation goals and $u$ is a function of three spatial variables. Many of the methods we shall present are not written as generalized for three dimensions. However, one may consider implementing the below schemes on tensor-product meshes, which is a straightforward formulation for some; one may also consider yet another splitting strategy in which separate schemes deal with only one-dimensional advection.

## 3.3 Flux Balance Methods

When we consider discretizing the full (6+1)-dimensional Vlasov equation (26), we first start by solving one-dimensional problems of the form

$$\frac{\partial f}{\partial t} + \frac{\partial (v f)}{\partial x} = 0. \qquad (29)$$

for some velocity $v(x,t)$, and $f(x,t)$, $x \in \mathbb{R}$. Although $v$ is independent of $x$ and $t$, we have assumed some generality by writing the equation in conservative form. We shall solve this system exactly using the method of characteristics. If we are given the solution at some time $t^n$, then we can easily write down the solution at time step $n+1$ as

$$f(x, t^{n+1}) = f(X(t^n; t^{n+1}, x), t^n). \qquad (30)$$

where we define $X(\,\cdot\,; t^{n+1}, x)$ as the solution to the ODE

$$\begin{aligned} \frac{\mathrm{d}X(s)}{\mathrm{d}s} &= v(X(s), s) \\ X(t^{n+1}) &= x. \end{aligned} \qquad (31)$$

Intuitively, $X(t^n; t^{n+1}, x)$ means 'the starting location at time $t^n$ of a particle that ends up at $x$ at time $t^{n+1}$ under the influence of the velocity field $v$'. Figure 5 furnishes a pictorial explanation.

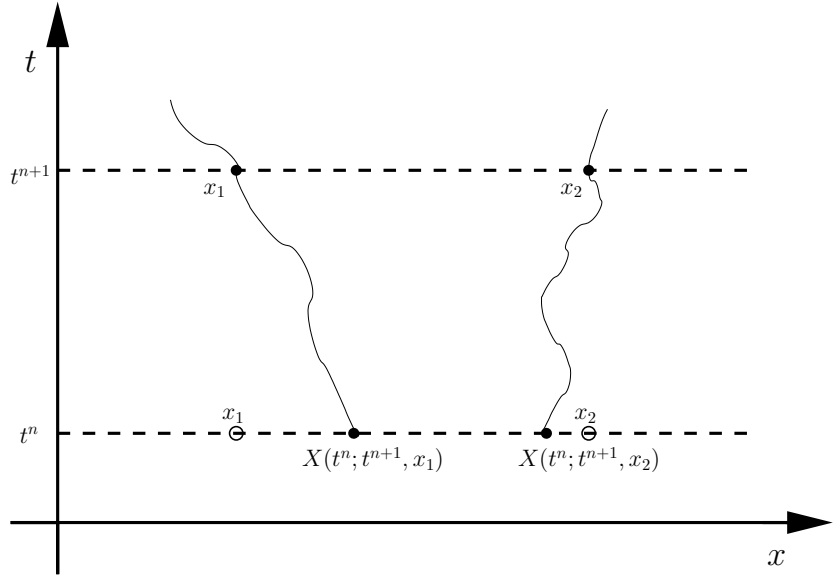

**Figure 5.** Pictorial representation of the characteristic flow governed by the evolution law (31).

Suppose that the coordinate direction $x$ is divided into a finite number of non-overlapping cells indexed by $i$ whose interfaces are given by the ordered nodal points $x_{i\pm 1/2}$, and whose centers are given by $x_i = \frac{1}{2}(x_{i-1/2} + x_{i+1/2})$. For the flux balance method (FBM), we shall consider the cell averages as degrees of freedom:

$$f_i^n = \frac{1}{\Delta x}\int_{x_{i-1/2}}^{x_{i+1/2}} f(t^n, x)\,\mathrm{d}x,$$

and introducing the flux terms

$$\Phi_{i+1/2}(t^n) = \int_{x_{i+1/2}-v\Delta t}^{x_{i+1/2}} f(t^n, x)\,\mathrm{d}x, \tag{32}$$

The flux $\Phi$ can be computed using the primitive of $f$:

$$\Phi_{i+1/2}(t^n) = F(t^n, x_{i+1/2}) - F(t^n, x_{i+1/2-v\Delta t}) \tag{33}$$

The primitive function $F$ satisfies the properties $F(t^n, x_{i+1/2}) - F(t^n, x_{i-1/2}) = \Delta x\, f_i^n$ and $F(t^n, x_{i+1/2}) = \Delta x \sum_{k=0}^{i} f_k^n =: w_i^n$. Assuming $v$ is divergence-free (which it is in the Vlasov equation), it follows from the conservation law formulation (29) that $f$ is conserved along characteristic boundaries. That is,

$$\int_{x_{i-1/2}}^{x_{i+1/2}} f(t^{n+1}, x)\,\mathrm{d}x = \int_{X(t^n; t^{n+1}, x_{i-1/2})}^{X(t^n, t^{n+1}, x_{i+1/2})} f(t^n, x)\,\mathrm{d}x, \tag{34}$$

which we may write in a finite volume framework as

$$f_i^{n+1} = f_i^n + \frac{\Phi_{i-1/2}(t^n) - \Phi_{i+1/2}(t^n)}{\Delta x}. \tag{35}$$

Equation (35) forms the basis for the solvers we shall discuss in this section. The flux terms $\Phi$ are computed using either (32) or (33), and the *modus operandi* for computing these fluxes is what differentiates each of the methods below.

### 3.3.1 The linear FBM

The FBM scheme (35) was introduced by Fijalkow in [14] where a linear interpolation procedure was used to determine a numerical reconstruction of the density $f_h$:

$$f_h(x) = f_i + (x - x_i) \frac{f_{i+1} - f_{i-1}}{2\Delta x}. \tag{36}$$

Using this reconstruction $f_h(x)$ at time level $t^n$, one can compute the fluxes explicitly using (32). The great advantage of this method is the ease of implementation. However, the main drawback is that positivity of the density is nowhere preserved, and spurious oscillations may develop.

### 3.3.2 The Positive and Flux Conservative Method

One well-known method to control positivity of the solution was introduced in [18]. The authors mention that using the property of the primitive function

$$F(t^n, x_{i+1/2}) = \Delta x \sum_{k=0}^{i} f_k^n, \tag{37}$$

one may use as large a centered finite volume stencil as desired to reconstruct an $n^{\text{th}}$ order primitive polynomial, which by differentiation yields an $n^{\text{th}}$ order accurate reconstruction of $f$, which may then be used to explicitly evaluate the integral in (32). However, there is again no way of preserving positivity with this method, and it will not even control spurious oscillations if the reconstructing stencil is not chosen in an adaptive way as in e.g. ENO [46].

In order to control the positivity, we adopt a third-order fixed-stencil reconstruction previously described, but nonlinear slope limiters $\varepsilon_i^\pm$ are introduced (see [17], [18] for details). The result is that the density reconstruction $f_h(x)$ has the form

$$\begin{aligned} f_h(x) &= f_i + \frac{\varepsilon_i^+}{6\Delta x^2} \big[ 2(x - x_i)(x - x_{i-3/2}) + (x - x_{i-1/2})(x - x_{i+1/2}) \big](f_{i+1} - f_i) \\ &\quad - \frac{\varepsilon_i^-}{6\Delta x^2} \big[ 2(x - x_i)(x - x_{i-3/2}) + (x - x_{i-1/2})(x - x_{i+1/2}) \big](f_i - f_{i-1}). \end{aligned} \tag{38}$$

The slope limiters are defined by

$$\varepsilon_i^\pm = \begin{cases} \min\left\{1, \dfrac{2 f_i}{f_{i\pm 1} - f_i}\right\} & \text{if } f_{i\pm 1} - f_i > 0, \\ \min\left\{1, -2\dfrac{f_\infty - f_i}{f_{i\pm 1} - f_i}\right\} & \text{if } f_{i\pm 1} - f_i > 0, \end{cases}$$

where $f_\infty = \max_i f_i$ is the maximum cell average value. Using this reconstruction for the density allows for integration to find the value of the primitive. The authors of [18] show that using this as the reconstruction of the density preserves both positivity and conservation of the average. That is, $f_h(x)$ satisfies

- *Conservation:* for every $i$, $\int_{x_{i-1/2}}^{x_{i+1/2}} f_h(t^n, x)\, \mathrm{d}x = \Delta x\, f_i^n$

- *Positivity:* for all $x$, $0 \leq f_h(t^n, x)$

With this well-behaved reconstruction for $f_h$ given by (38), one may then proceed to explicitly evaluate the integral (32) to advance the scheme (35); this method is called the positive and flux conservative (PFC) method. One can prove first-order convergence for this method. One-dimensional tests shown in [18] indicate that this scheme seems to control oscillations as well as ENO, and 2D tests confirm the validity of the procedure in higher dimensions.

A possible generalization of this method is to unify both the well-accepted oscillation-controlling properties of the ENO scheme with the provable positivity and conservation of the PFC method: the above method fixes a central stencil for the third-order reconstruction. One may consider instead using an adaptively-chosen stencil prescribed by the ENO procedure, and then use a reconstruction similar to (38) with the appropriately modified slope correctors. Although the PFC scheme has already heuristically shown its ability to control oscillations, this adaptation may prove a profitable union of the positivity, conservation, and oscillation-free properties.

### 3.3.3 PWENO Reconstruction

Very recently, researchers have applied the 'pointwise' WENO (PWENO) reconstruction procedure outlined in [44] to reconstruct the value of the primitive function to evaluate (33). They present their formulation and results in [10]. We shall briefly describe the PWENO procedure here.

We recall that in our finite volume framework, the cell boundaries are denoted by $x_{i-1/2}$ and our unknowns are the cell averages $f_i^n$. Equation (33) requires the evaluation of $F(t^n, x)$ for some value of $x$. Note that we can explicitly calculate $F(t^n, x_{i-1/2})$ for all $i$ using (37). Thus, our PWENO problem is: calculate $F(x)$ given $F(x_{i-1/2})$ for all $i$. Given a point $x \in (x_{i-1/2}, x_{i+1/2})$ we assume that we fix some 'ambient' stencil $S$: the collection of degrees of freedom which we will use to create a reconstruction. In this example, we shall use an ambient 5-point stencil. We associate with this ambient stencil three smaller 3-point stencils $S_r$ which are contained inside the ambient stencil $S$. See Figure 6.

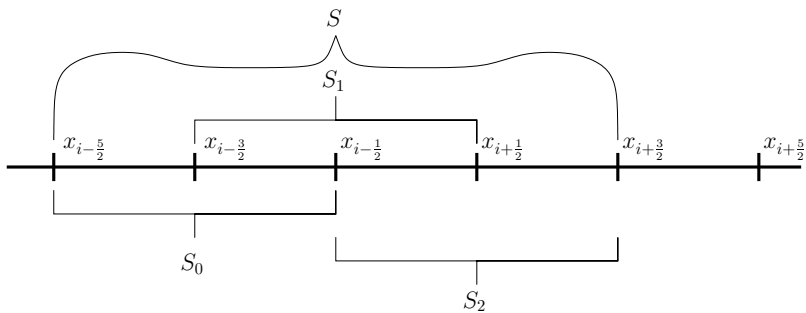

**Figure 6.** Depiction of a stencil setup for an ambient 5-point stencil comprised of three 3-point stencils.

The ENO/WENO schemes first introduced were concerned with reconstructing the function $f$ at the point $x = x_{i\pm 1/2}$. The algorithm requires the construction of polynomials approximating the primitive $F$ on each sub-stencil. We denote $P_r(x)$ as the polynomial reconstruction from the degrees of freedom $F(x_{i-1/2})$ on stencil $S_r$. Then, for example,

$$P_0(x) = \sum_{k=0}^{2} F[x_{i-5/2}, ..., x_{i+l-5/2}] \prod_{l=0}^{k-1} (x - x_{i-5/2+l}), \tag{39}$$

is the reconstruction on stencil $S_0$ where the Newton divided differences are calculated via the definition

$$F[x_i, ..., x_{i+p}] = \frac{F[x_{i+1}, ..., x_{i+p}] - F[x_i, ...x_{i+p-1}]}{x_{i+p} - x_i},$$

and the index $i$ may be shifted to e.g. $i - 5/2$ for the calculation of $P_0(x)$. As a base case, we define the divided difference $F[x_i] := F(x_i)$. To then obtain the reconstruction for our desired function of interest, $f$, we can differentiate (39) to obtain

$$p_0(x) = \sum_{k=1}^{2} F[x_{i-5/2}, ..., x_{i+l-5/2}] \sum_{l=0}^{k-1} \prod_{\substack{m=0 \\ m \neq l}}^{k-1} (x - x_{i-5/2+m}).$$

We now have a representation for our function of interest which depends only the divided differences of $F$ of the first order or higher. We can exploit this by noticing that the first-order divided difference $F$ is

$$F[x_{i-1/2}, x_{i+1/2}] = \frac{F(x_{i+1/2}) - F(x_{i-1/2})}{x_{i+1/2} - x_{i-1/2}} = f_i.$$

In other words, we can find reconstructions for $f$, e.g. $p_0$, without considering the primitive function at all.

This procedure is used by the ENO method to determine which stencil $S_r$ admits the least-oscillatory reconstruction $p_i(x)$, which it then uses to evaluate the reconstructed value. This is done by starting with a stencil size of 1, and increasing the stencil size by including the adjacent point which has the smallest divided difference (a mesaure of the oscillatory behavior of the reconstructed polynomial). More points are added until the intended stencil size is reached (for this example, 3 points). Once the stencil is computed, determine the reconstructing polynomial.

A generalization of the ENO procedure is the WENO procedure, which uses information from the entire ambient stencil $S$, rather than just one sub-stencil $S_r$. The ENO procedure has a sharp criterion for using a particular reconstruction: either one stencil is used entirely, or it is not used at all. The disadvantage of this is that, in the example in Figure 6, we use information from 5 points to reconstruct something with 3 points of accuracy. In an effort to increase accuracy for a given ambient stencil, WENO introduces the convex weights $d_r$ which are chosen so that

$$f_h(x_{i-1/2}) = \sum_{r=0}^{2} d_r p_r(x_{i-1/2}) + O(\Delta x^5),$$

and

$$\sum_{r=0}^{2} d_r = 1.$$

One can compute these weights readily based on Taylor expansions. However, we have simply computed a high-order reconstruction. This is what we would like the convex combination to behave like in smooth regions. To include some of the oscillation-reducing qualities of the ENO procedure in the vicinity of shocks, we introduce the smoothness indicators $\beta_r$ associated with each reconstruction $p_r(x)$, which we define as

$$\beta_r = \sum_{l=1}^{k-1} \int_{x_{i-1/2}}^{x_{i+1/2}} \Delta x^{2l-1} \left( \frac{\partial^l p_r(x)}{\partial x^l} \right)^2 dx,$$

where $k$ is the number of cells used in the reconstruction process. This indicator is a Sobolev-type seminorm based on the heuristic observation that highly oscillatory functions harbor high-magnitude derivatives. Thus, a higher value of $\beta_r$ indicates that a smaller contribution from $p_r(x)$ should be considered. Let $\varepsilon > 0$ be some small (e.g. $10^{-6}$) parameter which prevents machine overflow, and define

$$\tilde{\omega}_r = \frac{d_r}{(\varepsilon + \beta_r)^2}, \tag{40}$$

and

$$\omega_r = \frac{\tilde{\omega}_r}{\sum_{l=0}^{k-1} \tilde{\omega}_l}. \qquad (41)$$

In smooth regions the $\beta_r$ are all nearly equal so that the stencils are all weighed approximately equally. In oscillatory regions, oscillatory reconstructions have very small weights. The $\omega_r$ are the convexly normalized to preserve consistency. The WENO procedure uses

$$p(x_{i-1/2}) = \sum_{r=0}^{k-1} \omega_r p_r(x_{i-1/2}) \qquad (42)$$

as the reconstructed value. The WENO algorithm can be summarized as follows: we wish to determine $f(x_{i-1/2})$ using the point values of $F(x_{i-1/2})$:

1. Compute $\beta_r$ for each reconstruction $p_r$.

2. Calculate the $\omega_r$ for each reconstruction using (40) and (41).

3. Use (42) to evaluate $f(x_{i-1/2}) = p(x_{i-1/2})$.

This procedure has the weights tailored to admit high accuracy when the reconstruction is used at cell interfaces. Suppose that we wish to reconstruct function values which are not on a cell interface. This is where the 'pointwise' WENO (PWENO) procedure can be applied. The PWENO procedure is the same as the WENO procedure in spirit. However, some care must be taken since the desired reconstruction location $x$ is not necessarily located at the interface of any element. (The solution to (31) is not necessarily any grid point.) To mold the WENO calculation into the PWENO formulation, all that must be done is to augment the weights $d_r \rightarrow d_r(x)$ and the smoothness indicators $\beta_r \rightarrow \beta_r(x)$ so that they are functions of the location of reconstruction. The weights $d_r$ were derived from a high-order Taylor expansion, so solving for the $d_r$ is in principle possible, but it is not computationally trivial. The qualitative explanation of this is relatively straightforward, but the quantitative description and details of this method are quite complex, and we refer the reader to [44] and [10] for details.

The procedure for solving (29) is then clear: first we solve (31) for the value $X(t^n; t^{n+1}, x_{i-1/2})$. Then we determine an ambient stencil around $X(t^n; t^{n+1}, x_{i-1/2})$ and perform the PWENO process to reconstruct the value $F(t^n, X(t^n; t^{n+1}, x_{i-1/2}))$ given the pointwise primitive values defined by (37). The difference in this application is that we are interested in reconstructing the primitive of the original function $f$. However, the method for choosing the stencil remains the same.

Unfortunately, due to the relatively complicated nature of the PWENO reconstruction procedure, one cannot prove any of the qualities such as positivity or conservation for this type of flux-balance method. However, the oscillation-reducing properties of the WENO method are quite desirable for the Vlasov equation.

## 3.4 The Semi-Lagrangian Method

The semi-Lagrangian method uses the same motivation as is shown in Figure 5: the solution to the hyperbolic problem (29) can be determined by tracing the characteristics back in time using (31). The procedure is as follows: given $\{f(t^n, x_i)\}_{i \in I}$

- *Follow the characteristics* - Compute $\{X(t^n; t^{n+1}, x_i)\}_{i \in I}$ using (31).

- *Interpolate* - Use the known values $\{f(t^n, x_i)\}_{i \in I}$ to interpolate the values $\{f(t^n, X(t^n; t^{n+1}, x_i))\}_{i \in I}$, set $f(t^{n+1}, x_i) = f(t^n, X(t^n; t^{n+1}, x_i))$. See Figure 5.

The first step is relatively easy to accomplish, especially since the velocity field in (31) is constant for the $x$-advection. However, the second step of interpolation is quite a challenge. Indeed, this step is what differentiates each semi-Lagrangian method from the next. We present a few algorithms below and refer the reader to [47] and [17] for a more in-depth survey of these methods.

### 3.4.1 Lagrange Interpolations

The first semi-Lagrangian method was introduced by Cheng and Knorr [11]. Their attempts include linear, quadratic, and cubic interpolation. For a general $(2m+1)$th-order polynomial approximation of $f(t^n, \cdot)$ at $x$, we first determine $i$ such that $x \in [x_i, x_{i+1}]$. Then we fit a Lagrangian polynomial of order $2m+1$ with the nodal abscissae at $\{x_{i-m}, ..., x_{i-1}, x_i, ..., x_{i+m-1}\}$. This is given by the formula

$$f_h(x) = f(t^n, x_{i-m}) + \sum_{k=1}^{2m+1} f[x_{i-m}, ..., x_{i-m+k}] \prod_{l=0}^{k} (x - x_{i-m+1}),$$

where

$$f[x_i, ..., x_{i+p}] = \frac{1}{p!} \frac{f[x_{i+1}, ..., x_{i+p}] - f[x_i, ...x_{i+p-1}]}{x_{i+p} - x_i},$$

are the recursively-defined Newton divided differences, with $f[x_i] = f(t^n, x_i)$. This method relies on a centered stencil and is thus limited to odd-degree interpolation, but an even-degree interpolant may easily be implemented as well.

The authors in [11] observe that linear interpolation is too dissipative and maintain that cubic interpolation is satisfactory. Numerical simulations in [17] indicate that increasing interpolation orders (tests were performed up to a 9th degree interpolation) decreases the amount of dissipation and also decreases phase error. However, increasing the degree of the interpolant requires more computational effort, and widens the reconstructing stencil. This high-order reconstruction method does not inherently preserve positivity or mass.

### 3.4.2 Hermite Interpolations

The previous section introduced an interpolant which was an approximation to point values of the density $f$. The Hermite interpolation method introduced in [17] proposes to also use derivatives for reconstruction. Consider for simplicity of exposition a cubic reconstruction. The reconstruction $f_h(x)$ for $x \in [x_i, x_{i+1}]$ is sought so that

$$\begin{aligned} f_h(x_i) &= f(t^n, x_i) & f'_h(x_i) &= \frac{\partial f}{\partial x}(t^n, x_i) \\ f_h(x_{i+1}) &= f(t^n, x_{i+1}) & f'_h(x_{i+1}) &= \frac{\partial f}{\partial x}(t^n, x_{i+1}) \end{aligned}$$

The values $f(t^n, x_i)$ and $f(t^n, x_{i+1})$ are known, so these two constraints pose little problem. However, it is not straightforward how to calculate the derivative $\frac{\partial f}{\partial x}$ at the grid points. A method is introduced in [36] by which a second transport equation is introduced for which the unknowns are the derivatives $\frac{\partial f}{\partial x}$ evaluated at the grid points. This new transport equation reads

$$\frac{\partial (\partial_x f)}{\partial t} + \frac{\partial [v(\partial_x f)]}{\partial x} = 0,$$

where we have assumed that the velocity field $v$ is divergence-free. The disadvantage of this method is that one must solve another transport equation in tandem with the original, which increases computational cost and memory requirements, especially if we must do this in (6+1)-dimensional space.

An alternative is proposed by Filbet et al. in [17], where the authors approximate the derivative as

$$\frac{\partial f}{\partial x}(x_i) = \frac{1}{12 \Delta x} [8(f(x_{i+1}) - f(x_{i-1})) - (f(x_{i+2}) - f(x_{i-2}))], \tag{43}$$

which is a well-known fourth-order accurate finite-difference formula on an equidistant mesh. Use of this approximation renders the algorithm relatively easy to implement and makes it cost-effective. Numerical results show that this method performs on par with a 3rd or 5th order Lagrange interpolant.

Of course one may imagine a generalization of this Hermite interpolation method. Indeed, consider some integer $k \geq 0$ and some stencil size $s \geq 1$. Then given that $x \in [x_i, x_{i+s}]$, we could attempt to find a reconstruction $f_h(x)$ which satisfies

$$f_h^{(q)}(x_r) = f^{(q)}(x_r) \quad \forall \, q = 0, 1, ...k \text{ and } r = 0, 1, ...s. \tag{44}$$

Formulating this scheme is not too difficult to do, but implementing it is nontrivial. How do we compute $f^{(q)}(x_r)$? Finite-difference formulas of the form (43) are not readily available for any order derivative with any stencil. Additionally, these difference formulas may produce spurious oscillations when used. Indeed, implementing schemes of the form (44) seem at this time academic in interest only.

### 3.4.3 PWENO Reconstruction

We ask the reader to recall the PWENO reconstruction treated in section 3.3.3. This reconstruction forms the basis for this method. Given $x = X(t^n; t^{n+1}, x_i)$, suppose $x \in [x_j, x_{j+1}]$ for some $j$. We form the reconstruction $f_h(x) \approx f(t^n, X(t^n; t^{n+1}, x_i))$ from the point values $\{f(t^n, x_i)\}_{i \in I}$ with PWENO (see Section 3.3.3). This idea can be attributed to Carrillo and Vecil in [10]. Their numerical simulations show that using PWENO as a reconstruction method instead of the usual Lagrange interpolation for a one-dimensional advection problem produces much better results in the discrete $L^1$ norm, and also preserves the discrete total variation when the Lagrange interpolation method does not. More complicated two-dimensional phase space tests also indicate that the method performs very well when compared to other semi-Lagrangian or flux balance methods.

## 3.5 Other Methods

In the previous sections we have summarized some of the methods which have been developed to solve the Vlasov equation (26). So far these methods have been lumped into the general category of FBM solvers and semi-Lagrangian solvers. However, there are other methods that have been presented in the literature to solve this equation while maintaining positivity and conservation. Gamba and Proft [20] have applied a discontinuous Galerkin (DG) method to solve the problem, Arakawa [1] introduced a finite-difference method in 1966 for incompressible fluid flow which was later adapted by Filbet and Sonnendrücker in [17] for the Boltzmann equation, and Klimas and Farrell have explored spectral method representations in [33] and [32]. In this section, we shall briefly describe these three methods and their properties.

### 3.5.1 An Adapted Finite-Difference Scheme

In 1966, Arakawa [1] introduced a finite-difference method to solve the incompressible Euler equation of fluid flow in two dimensions. He showed that his scheme conserved the 'mean Jacobian', squared vorticity, and kinetic energy. Later in 2003, Filbet [17] adapted Arakawa's scheme for Vlasov-Poisson models and noted that the conserved quantities for the Euler equations which Arakawa derived corresponded to the quantities mass, energy, and mean-square of $f$, respectively. We follow Filbet's adaptation here: consider the $1+1$ phase-space Vlasov-Poisson model:

$$\frac{\partial f}{\partial t} + v \frac{\partial f}{\partial x} + \frac{\partial \phi}{\partial x} \frac{\partial f}{\partial v} = 0$$

$$\nabla_x^2 \phi = \int f \, \mathrm{d}v - 1 \tag{45}$$

Define the auxilliary variables

$$\psi = \phi - \frac{v^2}{2}, \qquad J(\psi, f) = \frac{\partial \psi}{\partial x}\frac{\partial f}{\partial v} - \frac{\partial \psi}{\partial v}\frac{\partial f}{\partial x}, \tag{46}$$

which we can use to rewrite (45) as

$$\frac{\partial f}{\partial t} + J(\psi, f) = 0.$$

Having made this transformation, Filbet notices that this is exactly the equation which Arakawa solved for fluid problems. One can show that Arakawa's conserved quantities can be translated into Boltzmann-type conserved quantities as:

$$\int J_h(\psi, f)\,\mathrm{d}x\,\mathrm{d}v = 0 \longrightarrow \int f_h(t)\,\mathrm{d}x\,\mathrm{d}v = \int f_h(0)\,\mathrm{d}x\,\mathrm{d}v,$$

$$\int \psi_h\, J_h(\psi, f)\,\mathrm{d}x\,\mathrm{d}v = 0 \longrightarrow \int \psi_h(t)\,f_h(t)\,\mathrm{d}x\,\mathrm{d}v = \int \psi_h(0)\,f_h(0)\,\mathrm{d}x\,\mathrm{d}v, \tag{47}$$

$$\int f_h\, J_h(\psi, f)\,\mathrm{d}x\,\mathrm{d}v = 0 \longrightarrow \int f_h^2(t)\,\mathrm{d}x\,\mathrm{d}v = \int f_h^2(0)\,\mathrm{d}x\,\mathrm{d}v,$$

where the integrals over the discrete representations refer to sums of the nodal values. These three conserved quantities correspond to mass, energy, and the $L^2$ norm of $f$. We introduce some standard notation for spatial finite-difference methods:

$$D_i^0(f_h)_{i,j} = (f_h)_{i+1,j} - (f_h)_{i-1,j}, \quad D_i^+(f_h)_{i,j} = (f_h)_{i+1,j} - (f_h)_{i,j}, \quad D_i^-(f_h)_{i,j} = (f_h)_{1,j} - (f_h)_{i-1,j},$$

$$D_{i,j}^{+;-}(f_h)_{i,j} = (f_h)_{i+1,j} - (f_h)_{i,j-1}, \qquad D_{i,j}^{-;+}(f_h)_{i,j} = (f_h)_{i,j+1} - (f_h)_{i-1,j}.$$

In order to define the spatial differences which we shall use to approximate $J_h$, Arakawa first defines some intermediate Jacobians:

$$\mathbb{J}_{i,j}^{++}(f_h, \psi_h) = \frac{1}{4h^2}\big[D_i^0(f_h)_{i,j}\, D_j^0(\psi_h)_{i,j} - D_j^0(f_h)_{i,j}\, D_i^0(\psi_h)_{i,j}\big],$$

$$\mathbb{J}_{i,j}^{+\times}(f_h, \psi_h) = \frac{1}{4h^2}\big[D_i^0((f_h)_{i,j}\, D_j^0(\psi_h)_{i,j}) - D_j^0((f_h)_{i,j}\, D_i^0(\psi_h)_{i,j})\big],$$

$$\mathbb{J}_{i,j}^{\times+}(f_h, \psi_h) = \frac{1}{4h^2}\Big[(f_h)_{i+1,j+1}(D_j^+ - D_i^+)(\psi_h)_{i,j} - (f_h)_{i-1,j-1}(D_j^+ - D_i^+)(\psi_h)_{i-1,j-1} -$$

$$(f_h)_{i-1,j+1} D_{i,j}^{-;-}(\psi_h)_{i,j+1} + (f_h)_{i+1,j-1} D_{i,j}^{-;-}(\psi_h)_{i+1,j}\Big],$$

$$\mathbb{J}_{i,j}^{\times+}(f_h, \psi_h) = \frac{1}{8h^2}\Big[D_{i,j}^{0,0}(f_h)_{i,j}\,(D_j^0(\psi_h)_{i-1,j} - D_i^0(\psi_h)_{i,j-1}) -$$

$$(D_j^0(f_h)_{i-1,j} - D_i^0(f_h)_{i,j-1})\, D_{i,j}^{0,0}(\psi_h)_{i,j}\Big].$$

Using these finite-difference approximations, he forms the following approximations of the term $J(f, \psi)$:

$$\mathbb{J}_{i,j} = \alpha\,\mathbb{J}_{i,j}^{++} + \beta\,\mathbb{J}_{i,j}^{+\times} + \gamma\,\mathbb{J}_{i,j}^{\times+} + \delta\,\mathbb{J}_{i,j}^{\times\times}. \tag{48}$$

Finally, he proves the following properties:

**Proposition 1.** *(Finite difference accuracy, conservation [1])*

Let (48) be the finite difference approximation to $J(f, \psi)$ given in (46). The following hold:

1. $\alpha = \beta = \gamma = \frac{1}{3}$ and $\delta = 0$ yields a second-order accurate scheme which conserves the properties in (47).

2. $\beta = \gamma = \delta = \frac{1}{3}$ and $\alpha = 0$ yields a second-order accurate scheme which conserves the properties in (47).

3. $\alpha = \frac{2}{3}$, $\beta = \gamma = \frac{1}{3}$, and $\delta = -\frac{1}{3}$ yields a fourth-order accurate scheme which conserves the properties in (47).

Although it is clear that this method has desirable properties for the Vlasov equation, it is unclear how to generalize this method to higher-dimensional phase space. Indeed, Arakawa's original method was developed explicitly for the two-dimensional system, and so generalizations to six dimensions is not *a priori* straightforward.

### 3.5.2 A Spectral Method

A number of studies have been published regarding spectral methods discretizations for the Vlasov equation. We shall use some of the work by Klimas [33] and [32] to give a brief exposition of the method. We start with the $1+1$ phase-space Vlasov-Poisson model

$$\frac{\partial f}{\partial t} + v \frac{\partial f}{\partial x} + \frac{\partial \phi}{\partial x} \frac{\partial f}{\partial v} = 0, \tag{49}$$

Klimas uses a Strang-splitting method to advect the Vlasov equation: we assume that there is periodicity in the $x$-direction and that $f$ has zero support outside $|v| \leq V_{\max}$ for some $V_{\max} \geq 0$ (see Section 4). The $x$-advection satisfies

$$\frac{\partial f}{\partial t} + v \frac{\partial f}{\partial x} = 0, \tag{50}$$

and assuming a Fourier Transform in $x$

$$f(x, v, t) = \sum_{n=-N}^{N} \hat{f}_n(v, t) \, e^{-inx},$$

and similarly in $v$:

$$f(x, v, t) = \sum_{n=-N}^{N} \tilde{f}_n(x, t) \, e^{-inv},$$

then we may write the exact advection step for a time-step $\Delta t$ as

$$f(x, v, t + \Delta t) = \sum_{n=-N}^{N} \hat{f}(v, t) \, e^{iv\Delta t} e^{-inx}. \tag{51}$$

Thus, we can evolve the split stages of (49) by simply phase-shifting the modes $\hat{f}(v, t)$. Similarly, we assume that we have computed the appropriate force $L$ satisfying (49), so that a $v$-advection step takes the form

$$f(x, v, t + \Delta t) \approx \sum_{n=-N}^{N} \tilde{f}(x, t) \, e^{iL\Delta t} e^{-inv}. \tag{52}$$

Thus, up to the approximation in (52) and the the Strang splitting error, we can solve the system (49) using successive steps of (51) and (52) and the FFT. This is an extremely attractive concept: a spectral treatment of the unknown and accurate temporal advancement. However, there are many problems with this method: firstly, we have assumed that $f$ is periodic in $x$ and has compact support in $v$. For many problems (e.g. an equilibrium Maxwellian density) this is not satisfied. Additionally, the support of $f$ in $v$ increases as time is advanced. Thus, eventually we shall develop a numerical solution $f$ which has very oscillatory behavior; i.e. the high modes are flooded. This is very easily shown: equation (50) has the solution $f(x,v,t) = f_0(x - vt, v)$ where $f_0(x,v)$ is the initial condition. Then we can calculate

$$\frac{\partial f}{\partial v} = \frac{\partial f_0}{\partial v} - t \frac{\partial f_0}{\partial x},$$

and the second term is clearly unbounded in time, i.e. the velocity dependence in $f$ has steep gradients, which translates into energy in the high modes of the $v$-space Fourier Transform. This is the so-called *filamentation* problem. Filamentation causes corruption of the numerical solution due to the fundamental nature of the spectral method. A simple solution would be to filter the resulting solution in velocity space. The problem with filtering is that it produces unpredictable damage to the accuracy of the solution.

For $\nabla_x \bar{L} = \int \bar{f} \, dv - 1$, Klimas [33] showed that a particular type of Gaussian filtering in $v$-space produces a solution which, instead of solving the system (49), solves a new system. Let us introduce the filter that we shall use over $v$-frequency space:

$$G(\omega) = v_0 \sqrt{\frac{2}{\pi}} e^{-2(\omega v_0)^2},$$

where $v_0 > 0$ is a parameter. Note that the limit of $v_0 \to 0$, then there is no filtering. Suppose that we consider the exact solution $f$ to (49) and we apply $G(v)$ to its $v$-Fourier Transform, and then transform back to $v$-space to obtain a solution $\bar{f}$. I.e, $\bar{f}$ is defined by

$$\bar{f}(x,v,t) = \frac{1}{2\pi(2N+1)} \sum_{n=-N}^{N} \left( \sum_{n=-N}^{N} G\left(\frac{2\pi n}{V_{\max}}\right) \tilde{f}_n(x,t) e^{-inv} \right) e^{inv}.$$

Klimas then shows that this new function (or, rather, its continuous Fourier Transform counterpart) satisfies the evolution system

$$\frac{\partial \bar{f}}{\partial t} + v \frac{\partial \bar{f}}{\partial x} + \bar{L} \frac{\partial \bar{f}}{\partial v} = -v_0^2 \frac{\partial^2 \bar{f}}{\partial x \, \partial v},$$

$$\frac{\partial}{\partial t}\left(\frac{\partial \phi}{\partial x}\right) = -\int v \bar{f} \, dx \, dv.$$

(53)

Klimas also showed that the parameter $v_0$ can be chosen dependent on the initial condition so that filamentation does not develop, i.e. so that one can obtain a provable bound on the velocity derivative of the function $\bar{f}$. In addition, $\bar{f}$ conserves the same mass and momentum as $f$:

$$\int f \, dx \, dv = \int \bar{f} \, dx \, dv,$$

$$\int f v \, dx \, dv = \int \bar{f} v \, dx \, dv.$$

The advantages of using the solution to the new system (53), or rather, of doing a specific kind of filtering to the original system (49), are provable prevention of velocity-space filamentation while preserving certain moments of the distribution function.

### 3.5.3 A Discontinuous Galerkin Method

Discontinuous Galerkin (DG) methods have recently been used in a variety of applications in numerical analysis [27]. DG methods have the advantage of combining the spectral-type convergence of modal representations with the finite-volume robustness for handling shock waves in hyperbolic systems. Very recently, Gamba and Proft have proposed a DG scheme specifically designed for linear Vlasov-Boltzmann models [20]. Their scheme is as follows:

Assume we have a space $\Omega \subset \mathbb{R}^3 \times \mathbb{R}^3$ representing a subspace of $(x, v)$ phase space with finite volume. We assume $\Omega$ has rectangular boundaries such that $x_i \in [0, L_i]$ and $v_i \in [-V_i, V_i]$ for $i = 1, 2, 3$. We assume a regular finite-element triangulation $\mathcal{T}_h = \{\omega\}$ of the space $\Omega$. $\mathcal{T}_h$ has boundary faces $\Gamma \stackrel{\circ}{=} \cup_{\omega \in \mathcal{T}_h} \partial \omega$, and we split $\Gamma$ up into seven regions defined by

$$
\begin{aligned}
&\Gamma_i && \text{The set of interior faces: } \Gamma \backslash \partial \Omega, \\
&\Gamma_{0^-} && \text{The set of faces } \gamma \in \Gamma \text{ representing } x_i = 0, v_i < 0, \\
&\Gamma_{0^+} && \text{The set of faces } \gamma \in \Gamma \text{ representing } x_i = 0, v_i > 0, \\
&\Gamma_{L^-} && \text{The set of faces } \gamma \in \Gamma \text{ representing } x_i = L_i, v < 0, \\
&\Gamma_{L^+} && \text{The set of faces } \gamma \in \Gamma \text{ representing } x_i = L_i, v > 0, \\
&\Gamma_{V^-} && \text{The set of faces } \gamma \in \Gamma \text{ representing } v_i = -V_i, \\
&\Gamma_{V^+} && \text{The set of faces } \gamma \in \Gamma \text{ representing } v_i = V_i.
\end{aligned}
$$

We define the space of discontinuous functions for some $k \geq 0$:

$$\Phi_h := \left\{ \phi \in L^2(\Omega) : \phi|_{(x,v) \in \omega} \in P^k(\omega) \, \forall \, \omega \in \mathcal{T}_h \right\},$$

where $P^k(\omega)$ represents the space of $k$'th order polynomials on the element $\omega$. We associate with each face of each element an outward-pointing normal vector $\hat{n}_\omega$. At the interface of an element $\omega$, the functions $\phi \in \Phi_h$ are discontinuous, so we specify the local value of the function as $\phi^-$ and the external value of the function as $\phi^+$. Similarly, $\hat{n}_\omega^-$ is the inward-pointing normal vector and $\hat{n}_\omega^+$ is the outward-pointing normal vector. Then for any element $\omega$ at any interface $\gamma \in \Gamma$ in the triangulation, and for any $\phi \in \Phi_h$, we define the average and jump operators

$$\{\phi\} = \frac{1}{2}(\phi^- + \phi^+), \qquad [\![\phi]\!] = \hat{n}_\omega^- \cdot \phi^- + \hat{n}_\omega^+ \cdot \phi^+.$$

Finally, we define a numerical upwind flux $\hat{\phi}$ as

$$\hat{\phi} = \{\phi\} + \frac{1}{2} [\![\phi]\!] \tag{54}$$

Now we can rewrite the Vlasov equation (26) as

$$\frac{\partial f}{\partial t} + \nabla \cdot (A f) = 0, \tag{55}$$

where

$$A = \begin{pmatrix} v \\ L \end{pmatrix},$$

and the $\nabla$ operator acts on the full $(x, v)$ phase space. Then we can form the usual DG scheme, seeking $f_h \in \Phi_h$ by multiplying (55) by a test function $\phi \in \Phi_h$ and enforcing $L^2$ orthogonality of the residual with an integration by parts:

$$\sum_{\omega \in \Omega} \int_\omega \frac{\partial f_h}{\partial t} \phi \, \mathrm{d}x \, \mathrm{d}v - \sum_{\omega \in \Omega} \int_\omega A f_h \cdot \nabla \phi \, \mathrm{d}x \, \mathrm{d}v + \sum_{\omega \in \Omega} \int_{\gamma = \partial \omega} \phi \, (A f_h)^* \cdot \hat{n}_\omega^- \, \mathrm{d}\gamma = 0, \tag{56}$$

where we have yet to define the numerical flux $(A\,f_h)^*$. Indeed, the choice of numerical flux is at the heart of each DG method. Recalling the upwind flux definition in (54), Gamba uses the following form for the fluxes:

$$(A\,f_h)^* = \begin{cases} \widehat{A\,f_h} & \text{on } \Gamma_i, \\ -A\,f_h & \text{on } \Gamma_{V^\pm}, \\ A\,f_h & \text{on } \Gamma_{0^-} \text{ and } \Gamma_{L^+}, \\ A\,f_0 & \text{on } \Gamma_{0^+}, \\ A\,f_L & \text{on } \Gamma_{L^-}, \end{cases} \qquad (57)$$

where we have introduced $f_0$ and $f_L$ as the boundary conditions on appropriate $x$-inflow boundaries. These choices of fluxes corresponds to upwinding in the interior of the domain, $x$-inflow and $x$-outflow at the appropriate $x$-boundaries, and homogeneous Dirichlet conditions at $v$-boundaries. Gamba and Proft also prove the following properties of the scheme (56)-(57):

**Proposition 2.** *(DG Stability, Convergence, Conservation, Positivity [20])*

*The scheme defined by (56)-(57) for the system (55) has the following properties:*

1. $L^1$ *stability:* $\|f_h(t=T)\|_{L^1(\Omega)} \leq \|f_h(t=0)\|_{L^1(\Omega)}$

2. *Convergence:* $\|f - f_h\|_{L^1(\Omega)} \leq t\,h^k(C_0 + C_1)$, *where $C_0$ depends on the mesh and higher-order Sobolev norms of $f$, and $C_1$ depends on the mesh and higher-order Sobolev norms of $f(t=0)$.*

3. *Conservation:*
$$\frac{\partial}{\partial t}\int_\Omega f_h\,\mathrm{d}x\,\mathrm{d}v = 0.$$

4. *Positivity: if $k=0$ (piecewise-constant approximation) and $f(t=0) \geq 0$, then the numerical solution to the semi-discrete form of (56) satisfies $f_h \geq 0$ for $t \geq 0$.*

Gamba and Proft have implemented the scheme on $1+1$-dimensional phase space with encouraging results.

## 3.6 Summary

In this section we have surveyed many methods commonly used in the literature to solve the collisionless Boltzmann equation. Recall that the main challenges in solving the Vlasov equation are necessary application to 6-dimensional phase space, conservation of 'mass', and positivity of the density function. Each of these separate challenges has been addressed in the literature, and we have discussed some of the methods above. We have compiled the results obtained into Table 1.

$$\begin{array}{rcl} h & - & \text{Chacteristic degree-of-freedom scale} \\ p & - & \text{Order of polynomial approximation} \\ \text{MassC} & - & \text{Conservation of mass} \\ \text{MomC} & - & \text{Conservation of momentum} \\ \text{EnC} & - & \text{Conservation of energy} \\ \text{Positivity} & - & \text{Preservation of density positivity} \\ \text{Order} & - & \text{Order of convergence} \\ \text{Multi-D} & - & \text{Number of (phase-space) dimensions for which method can be applied} \end{array}$$

| Method | Ref. | MassC | MomC | EnC | Positivity | Order | Multi-D[d] |
|---|---|---|---|---|---|---|---|
| FBM-Linear | [11] | ✓ | – | – | – | $(h)$ | (✓) |
| FBM-PFC | [18] | ✓ | – | – | ✓ | $h$ | (✓) |
| FBM-WENO | [10] | ✓ | – | – | – | $h^p$ | (✓) |
| SL-Lagrange | [17] | (✓)[a] | – | – | – | $h^p$ | 1D |
| SL-Hermite | [17] | (✓)[a] | – | – | – | $h^p$ | 1D |
| SL-WENO | [10] | – | – | – | – | $h^p$ | (✓) |
| FD-Arakawa | [1] | ✓ | – | ✓ | – | $h^2, h^4$ | 2D |
| FFT-Klimas | [32] | (✓)[b] | – | – | – | $(h^\infty)^{[b]}$ | (✓) |
| DG-Gamba | [20] | ✓ | – | ✓ | (✓)[c] | $h^p, h^\infty$ | ✓ |

**Table 1.** Summary of properties of Vlasov solvers.

Notes:

a) The semi-Lagrangian method for these interpolations is mass-conservative only in the case of constant-coefficient advection with the use of a central stencil.

b) The FFT method is a spectral method and thus is, in principle, spectrally accurate. However, the FFT method assumes periodicity of the density function in velocity space, which is almost never the case. In addition, the filamentation-free correction introduces filtering, which is not mass conservative.

c) The DG method only has proven positivity for piecewise-constant ($p = 0$) approximations. Due to the hp-adaptive nature of the discontinuous Galerkin method, the order of convergence can be both the order of the polynomial ($h$-refinement) and also spectral ($p$-refinement).

d) Dimensional adaptivity denoted by (✓) means that direct high-dimensional application may prove tricky but the method can be generalized via tensor-products, whereas ✓ means it is dimensionally-adaptive with an unstructured grid.

## 4 Collision Operator

In the following, we will discuss numerical techniques for the second part of the Boltzmann equation's right hand side, the collision term. It is helpful to recall the forms of the collision operator in equations (2), (3), and (4). Since we will likely be timestepping the collision operator and the advection term separately, we reduce (1) to time evolution using *only* the collision term, ignoring transport:

$$\frac{\partial f}{\partial t} = \frac{1}{k_n} Q(f, f). \tag{58}$$

Our first observation here is that the collision operator is local in $x$. Without the transport term, $f(x)$ and $f(x^*)$ for $x \neq x^*$ evolve completely independently of each other. Therefore, we may limit our considerations in this section to one particular point in space. This allows easy (spatial-subdivision-based) parallelization of algorithms for treating the Boltzmann equation.

We will begin this section by reviewing a little bit of the theory surrounding the Boltzmann collision operator, and highlight several (often drastic) changes that were made to the rather onerous collision integral with the goal of simplifying it. Next, we will discuss the discretization of velocity space and, armed with that knowledge, we will review a few representative methods in more detail. We will focus on deterministic methods, noting that stochastic methods for the collision operator [37, 30] are a fast, but noisy alternative for problems where high accuracy is not of the utmost importance. Finally, we will discuss a time stepping method that is particularly well-suited to the requirements of the Boltzmann collision operator.

## 4.1 Theory and Possible Simplificiations

Let $d$ be the dimension of velocity space. The Boltzmann collision operator conserves a few of its *moments*, namely mass, momentum and energy:

$$\int_{\mathbb{R}^d} Q(f,f) \begin{cases} 1 \\ v \\ |v|^2 \end{cases} \mathrm{d}v = 0.$$

The solution $f$ of (58) evolves towards a steady state called a *Maxwellian*:

$$f_\infty(v) := \frac{\rho}{(2\pi T)^{d/2}} \exp\left(-\frac{|V-v|^2}{2T}\right).$$

The Maxwellian depends on moments of the initial distribution:

**Density.**

$$\rho := \int_{\mathbb{R}^d} f(v)\mathrm{d}v$$

**Mean/Bulk Velocity.**

$$V := \frac{1}{\rho} \int_{\mathbb{R}^d} v f(v)\mathrm{d}v$$

**Temperature.**

$$T := \frac{1}{3\rho} \int_{\mathbb{R}^d} |V-v|^2 f(v)\mathrm{d}v$$

This evolution towards the Maxwellian holds even for the whole Boltzmann equation with the same Maxwellian $f_\infty$, but with slightly modified moments. This Maxwellian is then (sensibly) uniform in $x$.

Furthermore, the following technical theorem tells us about how (58) evolves toward equilibrium and is useful as a diagnostic for numerical methods:

**Theorem 3. (Boltzmann's $H$-Theorem)** *The collision operator $Q$ satisfies the inequality*

$$\int_{\mathbb{R}^d} Q(f,f)\log(f)\mathrm{d}v \leqslant 0.$$

The term

$$H(f) := \int_{\mathbb{R}^d} Q(f,f)\log(f)\mathrm{d}v$$

can be interpreted as negative entropy, and Theorem 3 implies that entropy can only increase.

### 4.1.1 The Bhatnagar-Gross-Krook Approximation

One approach to making the complicated Boltzmann collision operator more approachable is the *Bhatnagar-Gross-Krook ("BGK") Approximation* [49]. It rests on the fact that the effect of the collision term in (58) is mainly to push the distribution function towards an equilibrium state. Therefore, the idea here is to simply replace $Q^{\mathrm{lin}}(f,f)$ by the difference between $f$ and its equilibrium value, divided by a relaxation time controlling the speed of the evolution:

$$Q^{\mathrm{BGK}}(f) = \frac{1}{\tau}(f_\infty - f).$$

Here $\tau$ is the typical time scale associated with collision relaxation to the local equilibrium. In principle the relaxation time $\tau$ is a complicated functional of the distribution function $f$. One key simplification associated with BGK is the assumption of a constant value for this relaxation scale.

It is surprising that even the rather crude BGK approximation of the full Boltzmann equation is still powerful enough to derive the full set of "classical" hydrodynamic equations from it, including the continuity equation and the Navier-Stokes equation [49].

This implies that computations of the BGK-Boltzmann equation capture *at least* as much detail as classical CFD calculations. There exist classes of CFD methods making use of exactly this fact.

## 4.2 Discretizing Velocity Space

The first step in any numerical calculation is the choice of an appropriate discretization. While the consequences of choosing a spatial discretization are most directly felt when treating the transport term, the choice of discretization of velocity space has its major impact on the computation of the collision term. While many methods make do with only one velocity space representation, some actually use two different representations and constantly transform back and forth between them.

The largest problem in dealing with velocity space is the fact that it, in contrast to our computational resources, is infinite, and there is no 'preferred' or 'obvious' method of squeezing it into a finite representation. Some attempts [45] have been made to expand $f$ into functions of unbounded support (often Hermite polynomials), however none of these seem to have managed to move past the BGK-approximated collision term.

We will first try to find some justification for only allowing $f$ to have finite support. Let $\mathcal{B}(C, R)$ denote the ball of radius $R$ centered at $C$, and $\mathrm{supp}_v(f)$ denote the set $\{v: f(v) \neq 0\}$.

**Proposition 4.** *(Compact support of the density function, Lemma 1 in [16], proved as Proposition 2.1 in [39])*

*Suppose* $\mathrm{supp}_v(f) \subset \mathcal{B}(0, R)$. *Then*

1. $\mathrm{supp}_v(Q(f,f)(v)) \subset \mathcal{B}(0, R\sqrt{2})$.

2. *We define the change of variables* $g := v - v_*$. *If we suppose* $v \in \mathcal{B}(0, R)$, *then we do not change the value of $Q$ if we restrict* $g \in \mathcal{B}(0, 2R)$, *i.e.*

$$Q(f,f)(v) = \int_{\mathcal{B}(0,2R)} \int_{S^2} B(|g|, \cos\theta) \left[ f(v') f(v'_*) - f(v) f(v-g) \right] d\omega \, dg.$$

3. *With this limitation, we automatically also have* $v_* = v - g, v', v'_* \in \mathcal{B}(0, (2+\sqrt{2})R)$.

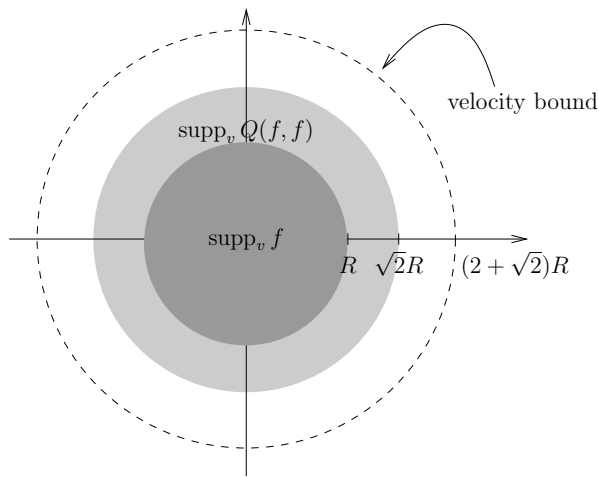

**Figure 7.** The Situation of Proposition 4.

This result does not encourage us to truncate the support of $f$: even if at timestep $n$ the support of $f$ falls within a ball of radius $R$, the support might grow by a factor of $\sqrt{2}$ with each evaluation of the collision operator.

### 4.2.1 Fourier Transforms in Velocity Space

Quite a few of the following methods use Cartesian Fourier transforms of the density function $f$ in velocity space. Since this is such a common theme, we discuss the implications of doing this beforehand. The aim in all cases is to use the fast Fourier transform (FFT) to speed up the computation of the expensive collision integral, similar to how convolution integrals are reduced from quadratic complexity to $O(N \log N)$ by the FFT.

We begin by discretizing velocity space $\mathbb{R}^3$ into the lattice $h_v \mathbb{Z}^3$:

$$v_k := h_v k, \quad h_v > 0, \quad k \in \mathbb{Z}^3. \tag{59}$$

Then, we approximate the integral

$$\varphi(\xi) = \int_{\mathbb{R}^3} f(v) e^{i(v \cdot \xi)} dv, \tag{60}$$

with the midpoint-rule like sum

$$\tilde{\varphi}(\xi) = h_v^3 \sum_{k \in \mathbb{Z}^3} f(v_k) e^{i(v \cdot \xi)}. \tag{61}$$

This immediately leads to the question at which points $\xi$ we would need to evaluate $\tilde{\varphi}(\xi)$. We resolve this by also picking a lattice in Fourier-velocity space:

$$\xi_j := h_\xi j, \quad h_\xi > 0, \quad j \in \mathbb{Z}^3. \tag{62}$$

In order to be able to carry out the transform (60) by using the FFT, we require

$$h_v h_\xi = \frac{2\pi}{n}, \quad n \in \mathbb{N}.$$

We see

$$\sum_{k \in \mathbb{Z}^3} f(v_k) e^{i\frac{2\pi}{n}((k + e_l n) \cdot j)} = \sum_{k \in \mathbb{Z}^3} f(v_k) e^{i\frac{2\pi}{n}(k \cdot j)} = \sum_{k \in \mathbb{Z}^3} f(v_k) e^{i\frac{2\pi}{n}(k \cdot (j + e_l n))},$$

where $e_l$ is the $l$th unit vector. We can therefore limit our grids (59) and (62) to $n^3$ elements each, and thereby automatically make our velocity space $2L$-periodic, with

$$L := \frac{n h_v}{2}.$$

The goal must be to pick $L$ as small as possible to improve resolution within $[-L, L]^3$, while picking it large enough to avoid spurious influence from periodic copies of the main mode of $f$ somewhere else in velocity space, an effect commonly known as *aliasing*.

Suppose $\mathrm{supp}_v(f) \subset \mathcal{B}(0, R_0)$. Taking into account Proposition 4, we know $\mathrm{supp}_v(Q(f, f)) \subset B(0, \sqrt{2} R_0)$, and that the relative velocity ($g$ above) is bounded in $B(0, 2R_0)$. Since we would really like our solutions within $\mathrm{supp}_v Q(f, f)$ to be correct, we need to make sure that these velocities keep a distance of $2R_0$ from the closest spurious non-zero location. Therefore, we obtain

$$L \geq \frac{1}{2}(\sqrt{2} R_0 + 2R_0 + R_0) = \frac{3 + \sqrt{2}}{2} R_0$$

as a guideline for the minimum period. Figure 8 aims to clarify this formula.

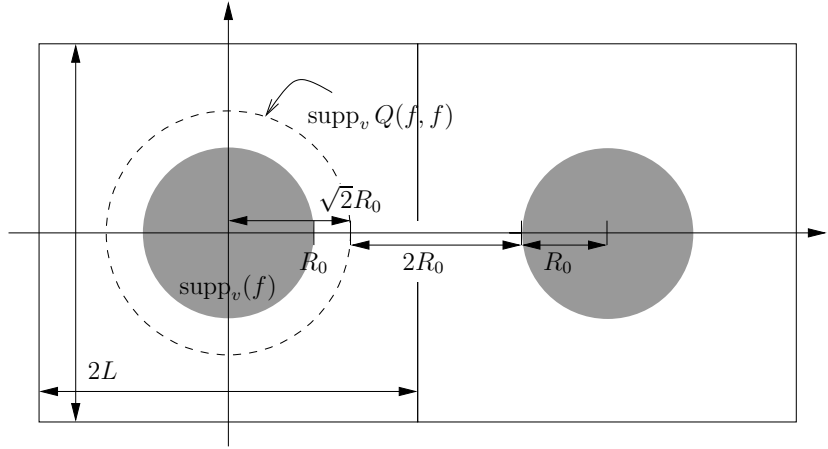

**Figure 8.** Considerations for mitigating the artificial $L$-periodicity in velocity space.

As a final note, some methods choose to change the definition (59) to originate at the (conserved) bulk velocity. This makes particular sense for situations where one is truly only interested in a solution of a space-homogeneous Boltzmann equation. (I.e. the velocity situation is assumed to be the same over all real space.) In this case, (59) is changed to

$$v_k := V + h_v k, \quad h_v > 0, \quad k \in \mathbb{Z}^3,$$

and (61) becomes

$$\tilde{\varphi}(\xi) = h_v^3 e^{i(V \cdot \xi)} \sum_{k \in \mathbb{Z}^3} f(v_k) e^{i(v \cdot \xi)},$$

with the bulk velocity $V$ defined as

$$\begin{aligned} V &:= \frac{1}{\rho} \int_{\mathbb{R}_v^3} v f(v,t) \mathrm{d}v = \frac{1}{\rho} \int_{\mathbb{R}_v^3} v f(v,0) \mathrm{d}v, \\ \rho &:= \int_{\mathbb{R}_v^3} f(v,t) = \int_{\mathbb{R}_v^3} f(v,0). \end{aligned}$$

### 4.3 A Fourier Method by Ibragimov and Rjasanow

The article [29] by I. Ibragimov and S. Rjasanow introduces a way of slicing the collision gain term (3) that is firmly rooted in the Fourier transform. To match the article and avoid confusion, we introduce the authors' notation for the Fourier transform of a function $z(v)$:

$$\begin{aligned} \mathcal{F}_{v \to \xi}[z(v)](\xi) &:= \int_{\mathbb{R}^3} z(v) e^{i(v \cdot \xi)} \mathrm{d}v, \\ \mathcal{F}_{\xi \to v}^{-1}[\hat{z}(\xi)](v) &:= \frac{1}{(2\pi)^3} \int_{\mathbb{R}^3} \hat{z}(\xi) e^{-i(v \cdot \xi)} \mathrm{d}\xi. \end{aligned}$$

Since the gain term is the most expensive part of approximating the collision operator (59), this method focuses on its treatment, beginning with the standard (substituted) form, recalling the definition of $g$ from Proposition 4.

$$\begin{aligned} Q_+(f,f)(v) &= \int_{\mathbb{R}^3} \int_{S^2} B(|g|, \cos\theta) \, f(v') \, f(v_*') \, \mathrm{d}\omega \, \mathrm{d}g \\ &= \int_{\mathbb{R}^3} \int_{S^2} B(|g|, \cos\theta) \, f\!\left(\tfrac{1}{2}(v + v_* + |g|\omega)\right) f\!\left(\tfrac{1}{2}(v + v_* - |g|\omega)\right) \mathrm{d}\omega \, \mathrm{d}g \\ &= \int_{\mathbb{R}^3} \int_{S^2} B(|g|, \cos\theta) \, f\!\left(\tfrac{1}{2}(v + v - g + |g|\omega)\right) f\!\left(\tfrac{1}{2}(v + v - g - |g|\omega)\right) \mathrm{d}\omega \, \mathrm{d}g \end{aligned}$$

$$= \int_{\mathbb{R}^3} \int_{S^2} B(|g|, \cos\theta) f\left(v + \frac{g}{2} + \frac{|g|\omega}{2}\right) f\left(v + \frac{g}{2} - \frac{|g|\omega}{2}\right) d\omega\, dg$$

$$= \int_{\mathbb{R}^3} \int_{S^2} B(|g|, \cos\theta) f\left(v - \frac{g}{2} + \frac{|g|\omega}{2}\right) f\left(v - \frac{g}{2} - \frac{|g|\omega}{2}\right) d\omega\, dg$$

We switch the integration of $g$ over $\mathbb{R}^3$ to spherical coordinates

$$Q_+(f,f)(v) = \int_0^\infty r^2 \int_{S^2} \int_{S^2} B(r, \omega \cdot \tilde{\omega}) f\left(v - \frac{r\tilde{\omega}}{2} + \frac{r\omega}{2}\right) f\left(v - \frac{r\tilde{\omega}}{2} - \frac{r\omega}{2}\right) d\omega\, d\tilde{\omega}\, dr,$$

revealing a very symmetric form. We substitute again, using $2u = r\omega$ (implying $r^2 dr\, d\omega = 8 du$), giving

$$Q_+(f,f)(v) = 8 \int_{\mathbb{R}^3} \int_{S^2} B(2|u|, \cos\theta) f(v - |u|\tilde{\omega} + u) f(v - |u|\tilde{\omega} - u)\, d\tilde{\omega}\, du.$$

While very similar to the expression we started with, we have unified the dependency of the $f$-terms on $\tilde{\omega}$. We omit the tilde over $\tilde{\omega}$:

$$Q_+(f,f)(v) = 8 \int_{\mathbb{R}^3} \int_{S^2} B(2|u|, \cos\theta) f(v - |u|\omega + u) f(v - |u|\omega - u)\, d\omega\, du.$$

The goal is now to extract this common dependency using the Fourier transform:

$$f(v - |u|\omega + u) f(v - |u|\omega - u) = \int_{\mathbb{R}^3} e^{iv\cdot(v - |u|\omega)} \mathcal{F}^{-1}_{z \to y}[f(z+u) f(z-u)] dy.$$

Altogether, we obtain

$$Q_+(f,f)(v) = 8 \int_{\mathbb{R}^3} \int_{S^2} B(2|u|, \cos\theta) \int_{\mathbb{R}^3} e^{iv\cdot(v - |u|\omega)} \mathcal{F}^{-1}_{z \to y}[f(z+u) f(z-u)] dy\, d\omega\, du.$$

Plucking this formula apart yields the representation

$$Q_+(f,f)(v) = \mathcal{F}_{y \to v}\left[\int_{\mathbb{R}^3} T(u,y) \mathcal{F}_{z \to y}[f(z-u) f(z+u)](u,y)\, du\right](v), \tag{63}$$

where the transformed kernel $T$ encapsulates the integration over the unit sphere, making it available for precomputation outside the timestepping loop:

$$T(u,y) = 8 \int_{S^2} B(2|u|, \cos\theta) e^{-i|u|y\cdot\omega} d\omega.$$

A few remarks are in order:

- It seems we have made a mess of the operator. Where we had integrals over five dimensions before, we now have *nine*. However, six of those are Fourier transforms that only contribute $O(N^3 \log N)$ complexity. Altogether, the final complexity of this method will turn out to be "only" $O(N^6 \log N)$.

- It is an important benefit that we are able to move the unwieldy evaluation of the unit sphere integral into a preprocessing step, where we can invest as much or as little time as we please.

- Some more analysis reveals that the kernel $T$ depends only on $|u|$, $|y|$ and $u \cdot y$. In the VHS case, it even only depends on $|u|$ and $|y|$.

- Carrying out the computations for (63) is reasonably straightforward, since we are only dealing with Cartesian integrals that need to be truncated appropriately. Large savings can be realized by taking advantage of the limited dependencies of $T$, if used properly. The reader is encouraged to look up the (somewhat messy) details in [29].

- This analysis only treats the (more expensive) gain term. Therefore, the loss term has to be computed separately. We refer to the source article [29] for some ideas.

## 4.4 The Integral Transform Method by Bobylev-Rjasanow

This method, introduced by Bobylev and Rjasanow in [6], is based on a clever reformulation of the collision operator. To start out, we once again use the substituted operator

$$Q(f,f) = \int_{\mathbb{R}^3} \int_{S^2} B(|v-v_*|, \cos\theta) \left[ f(v')f(v'_*) - f(v)f(v_*) \right] d\omega \, dv_*$$

$$(g := v - v_*) = \int_{\mathbb{R}^3} \int_{S^2} B(|g|, \cos\theta) \left[ f(v')f(v'_*) - f(v)f(v-g) \right] d\omega \, dg.$$

It turns out that this particular method works only for the hard spheres model ($B(|g|, \theta) = |g|$), but we will delay specializing for a bit, because our initial results will come in handy in Section 4.6.

### 4.4.1 Some Calculus Magic

**Lemma 5.** *(Lemma 1 in [6])*

*For any "nice enough" test function $\varphi : \mathbb{R}^3 \to \mathbb{R}$, we have*

$$\int_{S^2} \varphi(|g|\omega - g) d\omega = \frac{1}{|g|} \int_{\mathbb{R}^3} \delta\left(z \cdot g + \frac{1}{2}|z|^2\right) \varphi(z) dz.$$

We are replacing a two-dimensional (sphere) integral with a three-dimensional (volume) integral that includes a $\delta$-function.

**Lemma 6.** *(Part 1 of Lemma 2 in [6])*

*The collision operator can be represented as*

$$Q(f,f) = 4 \int_{\mathbb{R}^3} \int_{\mathbb{R}^3} \tilde{B}(y,z) \delta(z \cdot y) [f(v+z)(f(v+y) - f(v)f(v+y+z)] dy \, dz, \tag{64}$$

with

$$\tilde{B}(y,z) := \frac{1}{|y+z|} B\left(|y+z|, \frac{y \cdot (y+z)}{|y||y+z|}\right).$$

**Proof.** We define

$$\begin{aligned}
\varphi(|g|\omega - g) &:= f(v')f(v'_*) - f(v)f(v_*) \\
&= f\left(\frac{1}{2}(v+v_* + |g|\omega)\right) f\left(\frac{1}{2}(v+v_* - |g|\omega)\right) - f(v)f(v-g) \\
&= f\left(\frac{1}{2}(v+v-g+|g|\omega)\right) f\left(\frac{1}{2}(g+v_*+v_* - |g|\omega)\right) - f(v)f(v-g) \\
&= f\left(v + \frac{|g|\omega - g}{2}\right) f\left(v_* - \frac{|g|\omega - g}{2}\right) - f(v)f(v_*),
\end{aligned}$$

where $v$ and $v_*$ are parameters. Applying Lemma 5, we find

$$\int_{S^2} \varphi(|g|\omega - g) \mathrm{d}\omega = \frac{1}{|g|} \int_{\mathbb{R}^3} \delta\left(z \cdot g + \frac{1}{2}|z|^2\right) \varphi(z) \mathrm{d}z$$

$$\int_{S^2} [f(v')f(v'_*) - f(v)f(v_*)] \mathrm{d}\omega = \frac{1}{|g|} \int_{\mathbb{R}^3} \delta\left(z \cdot g + \frac{1}{2}|z|^2\right) \left[f\left(v + \frac{z}{2}\right) f\left(v_* - \frac{z}{2}\right) - f(v)f(v_*)\right] \mathrm{d}z.$$

Recognizing the left-hand side as part of the collision integral, we multiply the collision kernel $B$ and integrate to find

$$Q(f, f) = \int_{\mathbb{R}^3} \int_{\mathbb{R}^3} \frac{1}{|g|} B(|g|, \cos\theta) \delta\left(z \cdot g + \frac{1}{2}|z|^2\right) \left[f\left(v + \frac{z}{2}\right) f\left(v_* - \frac{z}{2}\right) - f(v)f(v_*)\right] \mathrm{d}z \mathrm{d}g.$$

Substituting $\tilde{z} := z/2$ (or $z = 2\tilde{z}$) and omitting the tilde sign right away gives us

$$\begin{aligned}
Q(f, f) &= 2^3 \int_{\mathbb{R}^3} \int_{\mathbb{R}^3} \frac{1}{|g|} B(|g|, \cos\theta) \delta\left(2z \cdot g + \frac{1}{2}|2z|^2\right) [f(v + z)f(v_* - z) - f(v)f(v_*)] \mathrm{d}z \mathrm{d}g \\
&= 8 \int_{\mathbb{R}^3} \int_{\mathbb{R}^3} \frac{1}{|g|} B(|g|, \cos\theta) \frac{\delta(z \cdot g + |z|^2)}{2} [f(v + z)f(v_* - z) - f(v)f(v_*)] \mathrm{d}z \mathrm{d}g \\
&= 4 \int_{\mathbb{R}^3} \int_{\mathbb{R}^3} \frac{1}{|g|} B(|g|, \cos\theta) \delta(z \cdot g + |z|^2) [f(v + z)f(v_* - z) - f(v)f(v_*)] \mathrm{d}z \mathrm{d}g,
\end{aligned}$$

taking into account $\delta(\alpha x) = \delta(x)/|\alpha|$. Next, we substitute $y := v_* - z - v = -z - g$ (or $v_* = y + z + v$, $g = -y - z$) and obtain

$$\begin{aligned}
Q(f, f) &= 4 \int_{\mathbb{R}^3} \int_{\mathbb{R}^3} \frac{1}{|g|} B(|g|, \cos\theta) \delta(z \cdot (z - y - z)) [f(v + z)f(y + z + v - z) - f(v)f(y + z + v)] \mathrm{d}z \mathrm{d}g \\
&= 4 \int_{\mathbb{R}^3} \int_{\mathbb{R}^3} \frac{1}{|g|} B(|g|, \cos\theta) \delta(z \cdot y) [f(v + z)f(v + y) - f(v)f(v + y + z)] \mathrm{d}z \mathrm{d}g,
\end{aligned}$$

proving the claim. $\square$

**Theorem 7.** *(Part 2 of Lemma 2 in [6])*

*The collision operator* for hard spheres *can also be represented as*

$$Q(f, f) = \int_{S^2} \int_{S^2} \delta(\omega_1 \cdot \omega_2) [\Phi(v, \omega_1) \Phi(v, \omega_2) - f(v) \Psi(v, \omega_1, \omega_2)] \mathrm{d}\omega_1 \mathrm{d}\omega_2, \tag{65}$$

where

$$\Phi(v, \omega) := \int_{-\infty}^{\infty} |\rho| f(v + \rho\omega) \mathrm{d}\rho, \tag{66}$$

and

$$\Psi(v, \omega_1, \omega_2) := \int_{-\infty}^{\infty} \int_{\infty}^{\infty} |\rho_1||\rho_2| f(v + \rho_1 \omega_1 + \rho_2 \omega_2) \mathrm{d}\rho_1 \mathrm{d}\rho_2. \tag{67}$$

**Proof.** Introduce spherical coordinates in (64) and use the identity

$$4\delta(\omega_1 \cdot \omega_2) = \delta(\omega_1 \cdot \omega_2) + \delta((-\omega_1) \cdot \omega_2) + \delta(\omega_1 \cdot (-\omega_2)) + \delta((-\omega_1) \cdot (-\omega_2)),$$

with the appropriate substitutions to extend the radial integration from $(0, \infty)$ to $(-\infty, \infty)$. $\square$

A few remarks are in order:

- First of all, the point is that (66) and (67) are integrals of convolution type and lend themselves well to computation using Fourier transform techniques.

- Next, a quick count reveals that we've grown an extra integration dimension when compared to (3) and (4) (six now, as compared to five before). This is however counteracted by the presence of the $\delta$-function in (65).

- A closer look at (65) shows that the integration is over all pairs of orthogonal unit vectors $\omega_1$, $\omega_2$. Recognizing this special structure should enable us to apply some cleverness in evaluating the integral.

- Expressions like (66) and (67) are not unheard of in mathematics. They resemble the X-ray and the Radon transform, respectively. [6] calls them *generalized X-ray* and *generalized Radon transforms*.

### 4.4.2 X-Ray to Fourier

To fix notation, we set the Fourier transform to be

$$\varphi(\xi) = \mathcal{F}[f](\xi) = \int_{\mathbb{R}^3} f(v) e^{i(v \cdot \xi)} \mathrm{d}v,$$

and its inverse

$$f(v) = \mathcal{F}^{-1}[\varphi](v) = \frac{1}{(2\pi)^3} \int_{\mathbb{R}^3} \varphi(\xi) e^{-i(v \cdot \xi)} \mathrm{d}\xi.$$

We begin by rewriting $\Phi$:

$$
\begin{aligned}
\Phi(v,\omega) &= \mathcal{F}^{-1}[\mathcal{F}[\Phi](\xi,\omega)](v,\omega) \\
&= \mathcal{F}^{-1}\left[\int_{\mathbb{R}^3} \Phi(v,\omega) e^{i(v \cdot \xi)} \mathrm{d}v\right](v,\omega) \\
&= \mathcal{F}^{-1}\left[\int_{\mathbb{R}^3} \int_{-\infty}^{\infty} |\rho| f(v + \rho\omega) \mathrm{d}\rho\, e^{i(v \cdot \xi)} \mathrm{d}v\right](v,\omega) \\
(\tilde{v} := v + \rho\omega) \quad &= \mathcal{F}^{-1}\left[\int_{\mathbb{R}^3} \int_{-\infty}^{\infty} |\rho| f(\tilde{v}) e^{i((\tilde{v} - \rho\omega) \cdot \xi)} \mathrm{d}\rho\, \mathrm{d}v\right](v,\omega) \\
&= \mathcal{F}^{-1}\left[\int_{\mathbb{R}^3} f(\tilde{v}) e^{i(\tilde{v} \cdot \xi)} \mathrm{d}v \int_{-\infty}^{\infty} |\rho| e^{-i(\rho\omega \cdot \xi)} \mathrm{d}\rho\right](v,\omega) \\
&= \mathcal{F}^{-1}[\varphi(\xi) d(\xi \cdot \omega)](v,\omega),
\end{aligned}
$$

with the definition

$$d(\zeta) := \int_{-\infty}^{\infty} |\rho| e^{-i\rho\zeta} \mathrm{d}\rho = \mathcal{F}[|\rho|](\zeta) = -\frac{2}{|\zeta|^2}.$$

Similarly for $\Psi$:

$$\Psi(v,\omega) = \mathcal{F}^{-1}[\varphi(\xi) d(\xi \cdot \omega_1) d(\xi \cdot \omega_2)](v, \omega_1, \omega_2).$$

The function $d$, however, is singular at $\zeta = 0$, and therefore not very suitable for our calculations. We therefore approximate

$$d(\zeta) \approx d_R(\zeta) := \int_{-R}^{R} |\rho| e^{-i\rho\zeta} \mathrm{d}\rho = 2\frac{\zeta R \sin(\zeta R) + \cos(\zeta R) - 1}{\zeta^2},$$

with $d_R(0) = R^2$. Note that this approximation has the net effect of truncating the integrals for $\Phi$ and $\Psi$ to the regions $[-R, R]$ and $[-R, R]^2$. We add the subscripts $Q_R$, $\Phi_R$ and $\Psi_R$ to denote that the approximate function $d_R$ is being used. The modified collision operator $Q_R$ still verifies several analytic properties of the original $Q$ such as the conservation laws, the $H$-theorem, and the equilibrium solutions.

This, of course, bears the question of how large $R$ should be chosen. We suppose that $\mathrm{supp}(f) \subset \mathcal{B}(0, R_0)$. Then, applying Proposition 4, we find that $|u| < 2R_0$, and $Q(f,f) = Q_R(f,f)$ for $R \geqslant R_0$. Therefore, if we choose a specific $R$, then the $L^\infty$ in $f$ error can be approximated by finding its maximum value of $f$ on the sphere $|v| = R/2$. Since $f$ tends to evolve into a Maxwellian $\exp(-\alpha |v|^2)$, we simply choose $R$ such that $\alpha R^2 \gg 4$, for the particular equilibrium in the problem to be solved.

### 4.4.3 Formula to Algorithm

All in all, we have so far analytically rewritten (59) into (65). What remains is to accurately approximate the various parts of (65). We proceed as follows:

1. Compute the Fourier transform $\varphi$ of $f$. The comments from Section 4.2.1 apply. In particular, [6] uses the grid shift by bulk velocity.

2. Carry out the integration
$$\int_{S^2} \int_{S^2} \delta(\omega_1 \cdot \omega_2) F(\omega_1, \omega_2) \mathrm{d}\omega_1 \, \mathrm{d}\omega_2.$$

   The main difficulty here is finding a suitable quadrature rule on the sphere. After that, the second integral collapses to a path integral over the circle $\{\omega : \omega \cdot \omega_1 = 0\} \subset S^2$, which is easy to carry out.

   At this stage, it is also useful to notice that the function
   $$F(\omega_1, \omega_2) := \Phi(v, \omega_1) \Phi(v, \omega_2) - f(v) \Psi(v, \omega_1, \omega_2)$$
   observes the symmetries
   $$F(\omega_1, \omega_2) = F(\pm \omega_1, \pm \omega_2) = F(\pm \omega_2, \pm \omega_1),$$
   which can be used for another substantial reduction of computational effort. For each value of $(\omega_1, \omega_2)$ found by the quadrature, do the following:

   a. Compute
   $$F_\Phi(\xi) := \varphi(\xi) d(\xi \cdot \omega), \qquad F_\Psi(\xi) := \varphi(\xi) d(\xi \cdot \omega_1) d(\xi \cdot \omega_2),$$
   from the (pre-computed) Fourier transform $\varphi$ of $f$.

   b. Compute the inverse Fourier transform of $F_\Phi$ and $F_\Psi$.

The method obtained in this way achieves second order convergence, and it conserves mass as-is. In order to conserve momentum and energy, it needs an ad-hoc modification to the axial terms of the Fourier transform that works by "making the smallest change" (in $L^2$ sense) that will conserve both quantities.

## 4.5 A Spectral Method by Pareschi, Russo, and Filbet

In [41], Pareschi et al. propose another Fourier-based scheme for the collision operator, and thus the remarks from Section 4.2.1 apply. For the purposes of their scheme, we assume that our function $f(v)$ has support on $\mathcal{B}(0, R_0)$, and we take the periodicity $L = \pi$ for simplicity. Naturally, we obtain the no-aliasing condition
$$\pi = L \geqslant \frac{3 + \sqrt{2}}{2} R_0,$$

as above. Let $K_N := \{-N, ..., N\}$. Not quite unlike 4.2.1, we expand $f$ in Fourier modes, with

$$f_N(v) := \sum_{k \in K_N^3} \hat{f}_k e^{i(k \cdot v)}, \qquad (68)$$

with the Fourier coefficients $\hat{f}_k$ given by

$$\hat{f}_k = \frac{1}{(2\pi)^3} \int_{[-\pi,\pi]^3} f(v) e^{-i(k \cdot v)} \mathrm{d}v.$$

### 4.5.1 Derivation

Our method is now defined by asserting that the residual of (58) is $L^2$-orthogonal to all Fourier modes:

$$\int_{[-\pi,\pi]^3} \left( \frac{\partial f_N}{\partial t} - \left[ Q^+(f_N, f_N) - L(f_N) f_N \right] \right) e^{-ik \cdot v} \mathrm{d}v \stackrel{!}{=} 0 \quad \text{for all } k \in K_N^3. \qquad (69)$$

Substituting (68) into (69) gives us (while setting $\mathcal{B} := \mathcal{B}(0, 2R_0)$)

$$\begin{aligned}
0 &\stackrel{!}{=} \int_{[-\pi,\pi]^3} \left( \frac{\partial f_N}{\partial t} - \left[ Q^+(f_N, f_N) - L(f_N) f_N \right] \right) e^{-i(k \cdot v)} \mathrm{d}v \\
&= \int_{[-\pi,\pi]^3} \left( \frac{\partial f_N}{\partial t} - \int_\mathcal{B} \int_{S^2} B(|v - v_*|, \cos\theta) \left[ f_N(v') f_N(v'_*) - f_N(v) f_N(v_*) \right] \mathrm{d}\omega \, \mathrm{d}v_* \right) e^{-i(k \cdot v)} \mathrm{d}v \\
&= \int_{[-\pi,\pi]^3} \left( \frac{\partial f_N}{\partial t} - \int_\mathcal{B} \int_{S^2} B(|g|, \cos\theta) \left[ f_N(v') f_N(v'_*) - f_N(v) f_N(v - g) \right] \mathrm{d}\omega \, \mathrm{d}g \right) e^{-i(k \cdot v)} \mathrm{d}v \\
&= \int_{[-\pi,\pi]^3} \left( \left( \sum_{k'' \in K_N^3} \frac{\partial \hat{f}_{k''}}{\partial t} e^{i(k'' \cdot v)} \right) - \int_\mathcal{B} \int_{S^2} B(|g|, \cos\theta) \left[ \left( \sum_{k' \in K_N^3} \hat{f}_{k'} e^{i(k' \cdot v')} \right) \left( \sum_{k'_* \in K_N^3} \hat{f}_{k'_*} e^{i(k'_* \cdot v'_*)} \right) - \right. \right. \\
&\qquad \left. \left. \left( \sum_{k'' \in K_N^3} \hat{f}_{k''} e^{i(k'' \cdot v)} \right) \left( \sum_{k_* \in K_N^3} \hat{f}_{k_*} e^{i(k_* \cdot (v-g))} \right) \right] \mathrm{d}\omega \, \mathrm{d}g \right) e^{-i(k \cdot v)} \mathrm{d}v \\
&= \int_{[-\pi,\pi]^3} \left( \left( \sum_{k'' \in K_N^3} \frac{\partial \hat{f}_{k''}}{\partial t} e^{i(k'' \cdot v - k \cdot v)} \right) - \int_\mathcal{B} \int_{S^2} B(|g|, \cos\theta) \left[ \sum_{k', k'_* \in K_N^3} \hat{f}_{k'} \hat{f}_{k'_*} e^{i(k' \cdot v' + k'_* \cdot v'_* - k \cdot v)} - \right. \right. \\
&\qquad \left. \left. \sum_{k'', k_* \in K_N^3} \hat{f}_{k''} \hat{f}_{k_*} e^{i(k'' \cdot v + k_* \cdot v - k_* \cdot g - k \cdot v)} \right] \mathrm{d}\omega \, \mathrm{d}g \right) \mathrm{d}v \\
&= \sum_{k'' \in K_N^3} \frac{\partial \hat{f}_{k''}}{\partial t} \int_{[-\pi,\pi]^3} e^{i(k'' \cdot v - k \cdot v)} \mathrm{d}v - \int_{[-\pi,\pi]^3} \int_\mathcal{B} \int_{S^2} B(|g|, \cos\theta) \left[ \sum_{k', k'_* \in K_N^3} \hat{f}_{k'} \hat{f}_{k'_*} e^{i(k' \cdot v' + k'_* \cdot v'_* - k \cdot v)} - \right. \\
&\qquad \left. \sum_{k'', k_* \in K_N^3} \hat{f}_{k''} \hat{f}_{k_*} e^{i(k'' \cdot v + k_* \cdot v - k_* \cdot g - k \cdot v)} \right] \mathrm{d}\omega \, \mathrm{d}g \, \mathrm{d}v \\
&= \sum_{k'' \in K_N^3} \frac{\partial \hat{f}_{k''}}{\partial t} \delta_{k'', k} - \int_{[-\pi,\pi]^3} \int_\mathcal{B} \int_{S^2} B(|g|, \cos\theta) \left[ \sum_{k', k'_* \in K_N^3} \hat{f}_{k'} \hat{f}_{k'_*} e^{i(k' \cdot v' + k'_* \cdot v'_* - k \cdot v)} - \right. \\
&\qquad \left. \sum_{k, k_* \in K_N^3} \hat{f}_{k''} \hat{f}_{k_*} e^{i(k'' \cdot v + k_* \cdot v - k_* \cdot g - k \cdot v)} \right] \mathrm{d}\omega \, \mathrm{d}g \, \mathrm{d}v \\
&= \frac{\partial \hat{f}_k}{\partial t} - \left[ \sum_{k', k'_* \in K_N^3} \hat{f}_{k'} \hat{f}_{k'_*} \int_\mathcal{B} \int_{S^2} \int_{[-\pi,\pi]^3} B(|g|, \cos\theta) \, e^{i(k' \cdot v' + k'_* \cdot v'_* - k \cdot v)} \mathrm{d}v \, \mathrm{d}\omega \, \mathrm{d}g - \right. \\
&\qquad \left. \sum_{k'', k_* \in K_N^3} \hat{f}_{k''} \hat{f}_{k_*} \int_\mathcal{B} \int_{S^2} \int_{[-\pi,\pi]^3} e^{i(k'' \cdot v + k_* \cdot v - k_* \cdot g - k \cdot v)} \mathrm{d}v \, \mathrm{d}\omega \, \mathrm{d}g \right]
\end{aligned}$$

$$
\begin{aligned}
&= \frac{\partial \hat{f}_k}{\partial t} - \Bigg[ \sum_{l,m \in K_N^3} \hat{f}_l \hat{f}_m \int_{\mathcal{B}} \int_{S^2} \int_{[-\pi,\pi]^3} B(|g|, \cos\theta)\, e^{i\left[l\cdot\left(\frac{1}{2}(v+v_*+|g|\omega)\right)+m\cdot\left(\frac{1}{2}(v+v_*-|g|\omega)\right)-k\cdot v\right]} \mathrm{d}v\, \mathrm{d}\omega\, \mathrm{d}g\ - \\
&\qquad \sum_{l,m \in K_N^3} \hat{f}_l \hat{f}_m \int \int_{S^2} \int_{[-\pi,\pi]^3} e^{i(l\cdot v+m\cdot v-m\cdot g-k\cdot v)} \mathrm{d}v\, \mathrm{d}\omega\, \mathrm{d}g \Bigg] \\
&= \frac{\partial \hat{f}_k}{\partial t} - \Bigg[ \sum_{l,m \in K_N^3} \hat{f}_l \hat{f}_m \int_{\mathcal{B}} \int_{S^2} \int_{[-\pi,\pi]^3} B(|g|,\theta)\, e^{i\left[l\cdot\left(\frac{1}{2}(v+v-g+|g|\omega)\right)+m\cdot\left(\frac{1}{2}(v+v-g-|g|\omega)\right)-k\cdot v\right]} \mathrm{d}v\, \mathrm{d}\omega\, \mathrm{d}g\ - \\
&\qquad \sum_{l,m \in K_N^3} \hat{f}_l \hat{f}_m \int_{\mathcal{B}} \int_{S^2} e^{i(-m\cdot g)} \delta_{l+m,k} \mathrm{d}\omega\, \mathrm{d}g \Bigg] \\
&= \frac{\partial \hat{f}_k}{\partial t} - \Bigg[ \sum_{l,m \in K_N^3,\, l+m=k} \hat{f}_l \hat{f}_m \int_{\mathcal{B}} \int_{S^2} B(|g|,\cos\theta)\, e^{i\left[l\cdot\left(\frac{1}{2}(-g+|g|\omega)\right)+m\cdot\left(\frac{1}{2}(-g-|g|\omega)\right)\right]} \mathrm{d}\omega\, \mathrm{d}g\ - \\
&\qquad \sum_{l,m \in K_N^3,\, l+m=k} \hat{f}_l \hat{f}_m \int_{\mathcal{B}} \int_{S^2} e^{i(-m\cdot g)} \mathrm{d}\omega\, \mathrm{d}g \Bigg] \\
&= \frac{\partial \hat{f}_k}{\partial t} - \Bigg[ \sum_{l,m \in K_N^3,\, l+m=k} \hat{f}_l \hat{f}_m \int_{\mathcal{B}} \int_{S^2} B(|g|, \cos\theta)\, e^{-ig\cdot \frac{l+m}{2} - i|g|\frac{m-l}{2}} \mathrm{d}\omega\, \mathrm{d}g\ - \\
&\qquad \sum_{l,m \in K_N^3,\, l+m=k} \hat{f}_l \hat{f}_m \int_{\mathcal{B}} \int_{S^2} e^{i(-m\cdot g)} \mathrm{d}\omega\, \mathrm{d}g \Bigg] \\
&= \frac{\partial \hat{f}_k}{\partial t} - \Bigg[ \sum_{l,m \in K_N^3,\, l+m=k} \hat{f}_l \hat{f}_m \hat{B}(l,m) - \sum_{l,m \in K_N^3,\, l+m=k} \hat{f}_l \hat{f}_m \hat{B}(m,m) \Bigg] \\
&= \frac{\partial \hat{f}_k}{\partial t} - \Bigg[ \sum_{m \in K_N^3} \hat{f}_{k-m} \hat{f}_m \hat{B}(k-m,m) - \sum_{m \in K_N^3} \hat{f}_{k-m} \hat{f}_m \hat{B}(m,m) \Bigg],
\end{aligned}
$$

where we have defined the *kernel modes* $\hat{B}(l,m)$ as

$$
\hat{B}(l,m) := \int_{\mathcal{B}(0,2R_0)} \int_{S^2} B(|g|, \cos\theta)\, e^{-ig\cdot \frac{(l+m)}{2} - i|g|\omega\cdot\frac{(m-l)}{2}} \mathrm{d}\omega\, \mathrm{d}g. \tag{70}
$$

We repeat the scheme for clarity:

$$
\frac{\partial \hat{f}_k}{\partial t} = \sum_{m \in K_N^3} \hat{f}_{k-m} \hat{f}_m \hat{B}(k-m,m) - \sum_{m \in K_N^3} \hat{f}_{k-m} \hat{f}_m \hat{B}(m,m). \tag{71}
$$

Again, now is the time to note down a few properties of this scheme:

- $\hat{B}(l,m)$ is a function of $|l-m|$, $|l+m|$, $(l-m)\cdot(l+m)$. (Proposition 3.2 in [41])

- For the VHS kernel, $\hat{B}(l,m)$ is a function of $|l-m|$, $|l+m|$. For this important case, there are only two parameters for $\hat{B}$, allowing it to be precomputed and stored in a two-dimensional array.

- Because we can pre-evaluate $\hat{B}$, evaluation of (71) takes $O(N^6)$ operations. In contrast, a straightforward method would take $O(N^6 N_a)$, where $N_a$ is the number of directions accounted for in a discretization of the unit sphere.

- The loss term is actually a convolution sum, and may therefore be evaluated in $O(N^3 \log N)$, making the gain term (sensibly) the largest consumer of computational resources.

It turns out that for the VHS kernels in two and three dimensions, the evaluation of (70) can be simplified drastically. We follow the two-dimensional derivation as an example.

### 4.5.2 Kernel Modes

Recall that for the VHS kernel, we are seeking a more explicit formula for the expression

$$\hat{B}_{\text{VHS}}(l,m) = C^\alpha \int_{\mathcal{B}(0,2R_0)} |g|^\alpha e^{-ig\cdot\frac{(l+m)}{2}} \underbrace{\int_{S^2} e^{-i|g|\omega\cdot\frac{(m-l)}{2}} d\omega}_{I_2(|g|,l-m):=} dg. \tag{72}$$

For the two-dimensional problem, we obtain

$$
\begin{aligned}
I_2(|g|,l-m) &= \int_{S^1} \exp\left(-i|g|\omega \cdot \frac{(m-l)}{2}\right) d\omega \\
&= \int_{S^1} \exp\left(-i|g|\omega \cdot e_1 \frac{|l-m|}{2}\right) d\omega \\
\left(r := |g|\frac{|l-m|}{2}, e_1 := (1,0)^T\right) &= \int_0^{2\pi} \exp(i r \cos\theta) d\theta \\
&= 2\int_0^{2\pi} \cos(r\cos\theta) d\theta \\
&= 2\pi J_0(r).
\end{aligned}
$$

Plugging this result into (72) and introducing polar coordinates yields

$$
\begin{aligned}
\hat{B}_{\text{VHS}}(l,m) &= 2\pi C^\alpha \int_{\mathcal{B}(0,2R_0)} |g|^\alpha \exp\left(-ig\cdot\frac{(l+m)}{2}\right) J_0\left(|g|\frac{|l-m|}{2}\right) dg \\
&= 2\pi C^\alpha \int_0^{2R_0} \rho^{1+\alpha} \int_0^{2\pi} \cos\left(|l+m|\frac{\rho}{2}\right)\cos\theta d\theta\, J_0\left(|g|\frac{|l-m|}{2}\right) d\rho \\
&= 4\pi^2 C^\alpha \int_0^{2R_0} \rho^{1+\alpha} J_0\left(|l+m|\frac{\rho}{2}\right) J_0\left(|g|\frac{|l-m|}{2}\right) d\rho \\
&= 4\pi^2 (2R_0)^{2+\alpha} \int_0^1 r^{1+\alpha} J_0(\xi r) J_0(\eta r) dr,
\end{aligned}
$$

with $\xi = |l+m|R_0$ and $\eta = |l-m|R_0$, allowing a very simple evaluation of $\hat{B}_{\text{VHS}}$. A similar analysis works in three dimensions.

### 4.5.3 Numerical Analysis

The question for the accuracy of the scheme is answered by the following:

**Theorem 8.** *(Theorem 5.3 in [41])*

Let $f \in L^2([-\pi,\pi]^3)$. Then

$$\left\|Q^{\text{tr}}(f,f) - \mathcal{P}_N Q^{\text{tr}}(f,f)\right\|_2 \leq C\left(\|f - f_N\|_2 + \frac{\|Q^{\text{tr}}(f_N, f_N)\|_{H^r_{\text{per}}}}{N^r}\right)$$

*for all $r \geq 0$, where $H^r_{\text{per}}$ is the Sobolev space of order $r$ containing only periodic functions, and $Q^{\text{tr}}$ represents the truncated collision operator.*

Stability can be proven once a certain smoothing (i.e. filtering, in Fourier space) operator is incorporated in the scheme. We refer to [41] for the details.

## 4.6 The Sub-$N^6$ Mouhot-Pareschi Method

In 2006, Mouhot and Pareschi introduced another method [19, 35] that cleverly combines the ideas from Sections 4.4 and 4.5. In doing so, it does achieve the important milestone of computing the entire collision operator in less than $O(N^6)$ operations, i.e. less than "any velocity coupled with any other".

Recall the statement of Lemma 6, slightly generalized for arbitrary dimensions:

**Lemma.** *The collision operator can be represented as*

$$Q(f,f) = \int_{\mathbb{R}^d} \int_{\mathbb{R}^d} \tilde{B}(y,z)\delta(z \cdot y)[f(v+z)f(v+y) - f(v)f(v+y+z)] \mathrm{d}y \, \mathrm{d}z, \tag{73}$$

with

$$\tilde{B}(y,z) := \frac{2^{d-1}}{|y+z|^{d-2}} B\left(|y+z|, \frac{y \cdot (y+z)}{|y||y+z|}\right).$$

The idea of this method is to use (73) instead of the "untreated" collision operator to derive a spectral method. Looking at the operator in this way yields better decoupling properties between the arguments of the operator, as we will see later.

The next job is to appropriately truncate the integration domains. Once again, we consider the periodic domain $[-L, L]^d$. If $\mathrm{supp}(f) \subset \mathcal{B}(0, R_0)$, then similar arguments as above show that we may truncate both integrations above to $\mathcal{B}(0, R)$, with $R = \sqrt{2} R_0$, while keeping $L \geqslant (3 + \sqrt{2})/2 \, R_0$ as the period. We obtain

$$Q^R(f,f) := \int_{\mathcal{B}(0,R)} \int_{\mathcal{B}(0,R)} \tilde{B}(y,z)\delta(z \cdot y)[f(v+z)f(v+y) - f(v)f(v+y+z)] \mathrm{d}y \, \mathrm{d}z.$$

Once truncated, we can derive the spectral scheme. The symbols $f_N$, $\hat{f}_k$, and $K_N$ are as in Section 4.5. We assume $L = \pi$. Again, we start by assuming that the residual of the evolution equation (58) is orthogonal to the space $\{v \mapsto e^{-ik \cdot v} : k \in K_N\}$.

$$\begin{aligned}
0 &\stackrel{!}{=} \int_{[-\pi,\pi]^d} \left(\frac{\partial f_N}{\partial t} - [Q^+(f_N, f_N) - L(f_N)f_N]\right) e^{-ik \cdot v} \mathrm{d}v \\
&= \int_{[-\pi,\pi]^d} \left(\sum_{k' \in K_N^d} \frac{\partial \hat{f}_{k'}}{\partial t} e^{i(k' \cdot v)} - Q^R(f_N, f_N)\right) e^{-ik \cdot v} \mathrm{d}v \\
&= \frac{\partial \hat{f}_k}{\partial t} - \int_{[-\pi,\pi]^d} [Q^R(f_N, f_N)] e^{-ik \cdot v} \mathrm{d}v \\
&= \frac{\partial \hat{f}_k}{\partial t} - \int_{[-\pi,\pi]^d} \int_{\mathcal{B}(0,R)} \int_{\mathcal{B}(0,R)} \tilde{B}(y,z)\delta(z \cdot y)[f_N(v+z)f_N(v+y) - f_N(v)f_N(v+y+z)] e^{-ik \cdot v} \mathrm{d}y \, \mathrm{d}z \, \mathrm{d}v \\
&= \frac{\partial \hat{f}_k}{\partial t} - \int_{\mathcal{B}(0,R)} \int_{\mathcal{B}(0,R)} \tilde{B}(y,z)\delta(z \cdot y) \int_{[-\pi,\pi]^d} \sum_{l,m \in K_N^3} \hat{f}_l \hat{f}_m \, e^{i[l \cdot (v+z) + m \cdot (v+y) - k \cdot v]} \\
&\qquad - \sum_{l,m \in K_N^3} \hat{f}_l \hat{f}_m e^{i[l \cdot v + m \cdot (v+y+z) - k \cdot v]} \mathrm{d}v \, \mathrm{d}y \, \mathrm{d}z \\
&= \frac{\partial \hat{f}_k}{\partial t} - \int_{\mathcal{B}(0,R)} \int_{\mathcal{B}(0,R)} \tilde{B}(y,z)\delta(z \cdot y) \sum_{l,m \in K_N^3} \hat{f}_l \hat{f}_m \delta_{l+m,k} \left[e^{i[l \cdot z + m \cdot y]} - e^{im \cdot (y+z)}\right] \mathrm{d}y \, \mathrm{d}z \\
&= \frac{\partial \hat{f}_k}{\partial t} - \sum_{l,m \in K_N^3, l+m=k} \hat{f}_l \hat{f}_m \int_{\mathcal{B}(0,R)} \int_{\mathcal{B}(0,R)} \tilde{B}(y,z)\delta(z \cdot y) \left[e^{i[l \cdot z + m \cdot y]} - e^{im \cdot (y+z)}\right] \mathrm{d}y \, \mathrm{d}z \\
&= \frac{\partial \hat{f}_k}{\partial t} - \sum_{l,m \in K_N^3, l+m=k} \hat{f}_l \hat{f}_m \hat{\beta}(l,m),
\end{aligned}$$

where
$$\hat{\beta}(l,m) := \int_{\mathcal{B}(0,R)} \int_{\mathcal{B}(0,R)} \tilde{B}(y,z)\delta(z\cdot y)\Big[e^{i[l\cdot z+m\cdot y]} - e^{im\cdot(y+z)}\Big]\mathrm{d}y\,\mathrm{d}z.$$

We repeat the scheme for clarity:
$$\frac{\partial \hat{f}_k}{\partial t} = \sum_{l,m\in K_N^3, l+m=k} \hat{f}_l \hat{f}_m \hat{\beta}(l,m), \tag{74}$$

and notice that the $\hat{\beta}(l,m)$ can be simplified further by writing
$$\hat{\beta}(l,m) = \beta(l,m) - \beta(m,m),$$

with
$$\beta(l,m) := \int_{\mathcal{B}(0,R)} \int_{\mathcal{B}(0,R)} \tilde{B}(y,z)\delta(z\cdot y)\Big[e^{i[l\cdot y+m\cdot z]}\Big]\mathrm{d}y\,\mathrm{d}z = \beta(m,l). \tag{75}$$

Recall that the scheme (71) that we derived in Section 4.5 had exactly the same shape. So how do we expect to save computational time here with a scheme that looks exactly the same? In Section 4.5 we truncated the integration in $g = y + z$, so our truncated operator there would look like this in our current representation:
$$Q_{\mathrm{old}}^R(f,f) = \int_{\mathcal{B}(0,R)} \int_{\mathcal{B}(0,R)} \underbrace{\mathbf{1}_{|y+z|\leqslant R}}_{\text{extra term}} \tilde{B}(y,z)\delta(z\cdot y)[f(v+z)f(v+y) - f(v)f(v+y+z)]\mathrm{d}y\,\mathrm{d}z.$$

It is exactly this extra characteristic function term that causes the variables $y$ and $z$ to couple unnecessarily, and this is exactly where we save (a lot of) time.

These savings obviously have to come from a careful analysis of the kernel modes $\beta(l,m)$ to allow simplification of (74). And that is exactly what we will tackle next.

We introduce spherical coordinates:
$$\beta(l,m) = \int_{S^{d-1}} \int_{S^{d-1}} \delta(\omega\cdot\omega') \int_0^R \int_0^R |\rho|^{d-2}|\rho'|^{d-2}\tilde{B}(|\rho|,|\rho'|)\Big[e^{i[l\cdot\rho\omega+m\cdot\rho'\omega']}\Big]\mathrm{d}\rho\,\mathrm{d}\rho'\,\mathrm{d}\omega\,\mathrm{d}\omega',$$

where we note that

a) due to the orthogonality enforced by the $\delta$-function, $\tilde{B}$ only depends on the moduli $\rho$ and $\rho'$,

b) the peculiarity with the $\delta$-function encountered in the proof of Lemma 5 reduces the exponent from the usual $\rho^{d-1}$ to $\rho^{d-2}$.

Now we extend the inner integrals from $(0,R)$ to $(-R,R)$, dividing by four to take care of the three extra terms introduced.

$$\beta(l,m) = \frac{1}{4}\int_{S^{d-1}} \int_{S^{d-1}} \delta(\omega\cdot\omega') \int_{-R}^R \int_{-R}^R |\rho|^{d-2}|\rho'|^{d-2}\tilde{B}(|\rho|,|\rho'|)\Big[e^{i[l\cdot\rho\omega+m\cdot\rho'\omega']}\Big]\mathrm{d}\rho\,\mathrm{d}\rho'\,\mathrm{d}\omega\,\mathrm{d}\omega'. \tag{76}$$

The last key ingredient in the efficient evaluation of (76) plugged into (74) is the "*decoupling assumption*"
$$\tilde{B}(y,z) = a(|y|)b(|z|) \quad \text{for} \quad y \perp z,$$

which, for us, boils down to
$$\tilde{B}(\rho,\rho') = a(\rho)b(\rho').$$

Obviously, Maxwellian molecules in 2D and Hard Spheres in 3D satisfy this assumption, because for these $\tilde{B}(y, z)$ is just a constant. Other cases where this can be made to happen are enumerated in the appendix of [35]. We specialize to three dimensions and find

$$\begin{aligned}\beta(l,m) &= \frac{1}{4}\int_{S^2}\int_{S^2}\delta(\omega\cdot\omega')\int_{-R}^{R}\rho a(\rho)e^{il\cdot\rho\omega}\mathrm{d}\rho\int_{-R}^{R}\rho' b(\rho')e^{im\cdot\rho'\omega'}\mathrm{d}\rho'\mathrm{d}\rho'\mathrm{d}\omega\mathrm{d}\omega' \\ &= \frac{1}{4}\int_{S^2}\int_{-R}^{R}\rho a(\rho)e^{il\cdot\rho\omega}\mathrm{d}\rho\int_{S^2\cap\omega^\perp}\int_{-R}^{R}\rho' b(\rho')e^{im\cdot\rho'\omega'}\mathrm{d}\rho'\mathrm{d}\omega'\mathrm{d}\rho\mathrm{d}\omega \\ &= \frac{1}{4}\int_{S^2}\varphi_a(l\cdot\omega)\int_{S^2\cap\omega^\perp}\varphi_b(m\cdot\omega')\mathrm{d}\omega'\mathrm{d}\rho\mathrm{d}\omega,\end{aligned}$$

where

$$\varphi_f(s) := \int_{-R}^{R}|\rho|f(|\rho|)e^{i\rho s}\mathrm{d}\rho,$$

and we notice $\varphi_f(s) = \varphi_f(-s)$. Therefore, it is sufficient to integrate over a hemisphere $(S^2)^+$, which gives us

$$\beta(l,m) = \int_{(S^2)^+}\varphi_a(l\cdot\omega)\int_{S^2\cap\omega^\perp}\varphi_b(m\cdot\omega')\mathrm{d}\omega'\mathrm{d}\rho\mathrm{d}\omega.$$

As a last step, we notice that the function

$$\omega \mapsto \varphi_a(l\cdot\omega)\int_{S^2\cap\omega^\perp}\varphi_b(m\cdot\omega')$$

is periodic on $(S^2)^+$. Therefore, even the simple rectangle rule quadrature with a parametrization

$$\theta := p\pi/M, \quad \varphi := q\pi/M,$$

of the hemisphere has infinite order:

$$\beta(l,m) \approx \tilde{\beta}(l,m) := \frac{\pi^2}{M^2}\sum_{p,q=0}^{M}\underbrace{\varphi_a(l\cdot\omega_{\theta,\varphi})}_{\alpha_{p,q}(l):=}\underbrace{\int_{S^2\cap\omega^\perp_{\theta,\varphi}}\varphi_b(m\cdot\omega')}_{\alpha'_{p,q}(m):=}.$$

The specific structure of $\beta(l,m)$ is the key to the speed of the method. Observe:

$$\begin{aligned}\frac{\partial\hat{f}_k}{\partial t} &= \frac{\pi^2}{M^2}\sum_{l,m\in K_N^3, l+m=k}\hat{f}_l\hat{f}_m\sum_{p,q=0}^{M}\alpha_{p,q}(l)\alpha'_{p,q}(m) \quad - \quad \text{loss} \\ &= \frac{\pi^2}{M^2}\sum_{m\in K_N^3}\hat{f}_{k-m}\hat{f}_m\sum_{p,q=0}^{M}\alpha_{p,q}(k-m)\alpha'_{p,q}(m) \quad - \quad \text{loss} \\ &= \frac{\pi^2}{M^2}\sum_{p,q=0}^{M}\sum_{m\in K_N^3}\hat{f}_{k-m}\alpha_{p,q}(k-m)\hat{f}_m\alpha'_{p,q}(m) \quad - \quad \text{loss}.\end{aligned}$$

Now this is a sum of $M^2$ convolutions, which can be carried out in $M^2N^3\log N$ operations, using an FFT, giving, to the best of our knowledge, the currently smallest asymptotic complexity of any deterministic method for the collision operator.

We close with a few remarks:

- As shown in [35], the method is spectrally accurate.

- The only "sophisticated" numerical technique required in the implementation of this method is the FFT. Therefore, the the net implementation effort is still reasonably limited.

- To be able to prove the usual conservation properties, the effective quadrature $(\mathcal{A}, d\mathcal{A}(\omega, \omega'))$ on the set

$$\{(\omega, \omega') \in S^{d-1} \times S^{d-1} : \omega \cdot \omega' = 0\}$$

  needs to be even, i.e. satisfy

$$(\omega, \omega') \in \mathcal{A} \quad \Rightarrow \quad (\omega, -\omega'), (-\omega, \omega'), (-\omega, -\omega') \in \mathcal{A},$$

  and the discrete measure $d\mathcal{A}(\omega, \omega')$ also needs to be even.

## 4.7 Summary

| Method | Ref. | Maxw. | HS | VHS | Other | Order | Complexity | C.Mass | C.Momt. | C.Energy |
|---|---|---|---|---|---|---|---|---|---|---|
| DSMC | [5] | ✓ | ✓ | ✓ | ✓ | "$N^{-1/2}$" | $N^3$ | ✓ | ✓ | ✓ |
| Ibragimov-R. | [29] | ✓ | ✓ | ✓ | – | $N^{-2}$ | $N^6$ | ✓ | (✓) [1] | (✓) [1] |
| general case [3] | [29] | – | – | – | ✓ | $N^{-2}$ | $N^6 \log N$ | ✓ | (✓) [1] | (✓) [1] |
| Bobylev-R. | [6] | – | 3D | – | – | $N^{-2}$ | $N^6 \log N$ | ✓ | (✓) [1] | (✓) [1] |
| Pareschi-Russo | [41] | ✓ | ✓ | ✓ | – | spectral | $N^6$ | ✓ | (✓) [2] | (✓) [2] |
| Mouhot-P. | [35] | 2D | 3D | – | – | spectral | $M^2 N^3 \log N$ | ✓ | (✓) [2] | (✓) [2] |

Legend:

**$N$.** $N$ is the number of velocity discretization points in one dimension.

**$M$.** $M$ is a fixed number of discretization angles.

**Order.** Implicitly in $O(\cdot)$-notation.

**Complexity.** Implicitly in $O(\cdot)$-notation.

> The collision operator is local in space, i.e. its value at $x$ only depends values at $x$. "Complexity" here is the number of operations required to calculate the collision operator at one point in space. To simulate the whole Boltzmann equation, one would need to do this $P^3$ times, where $P$ is the number of points in the one-dimensional spatial discretization.
>
> If we assume $P = M = N$, the fastest of these methods bring the naïve $N^{11}$ complexity down to $N^8 \log N$. Even if we assume only a measly 10 points per dimension, that is a win by almost three orders of magnitude!

**[1].** These conservations only hold if an extra correction is applied.

**[2].** Conservation holds up to the order of the scheme.

**[3].** The Ibragimov-Rjasanow scheme is split into two rows: The (faster) VHS case and the (slower) general case.

## 4.8 Time Discretization for the Collision Operator

Standard time discretizations, such as multi-stage Runge-Kutta or multi-step methods, while sometimes used to advance (58) in time, often fail to keep $f$ positive, while having the added drawback of restricting the timestep to rather small values.

Therefore, one is often better served taking advantage of the considerable freedom in choosing time discretizations afforded by the shape of (58). In particular, since each spatial point $x$ is independent of its surroundings, we are free to choose separate timesteps or even separate timestepping *methods* at each point in space.

Among just the standard methods, Filbet and Russo [16] favor multi-step over multi-stage methods for improved efficiency. In general however, they advocate the use of a class of special-purpose schemes, called a *time-relaxed (TR) schemes* (see [43]). These schemes are designed to maintain the positivity of $f$ while permitting large timesteps.

To construct a TR scheme, we begin by assuming that we want to advance in time an equation of the shape

$$\frac{\partial f}{\partial t} = \frac{1}{k_n}[P(f,f) - \mu f], \tag{77}$$

where $P(f,f)$ is a positive bilinear operator, and $\mu > 0$ is a constant. The point of this representation is that the constant $\mu$ puts an upper bound on the "negative part" (i.e. the loss term) of the collision operator.

For a Maxwellian gas, the kernel $B$ does not depend on velocity and therefore $Q$ already has the form (77). For more general kernels, we would need to find a $\mu$ to satisfy

$$\mu \geqslant L[f](v) \quad \text{for all } v \in \mathbb{R}^3. \tag{78}$$

Then we could rewrite (58) in the form (77) by setting

$$P(f,f) = Q^+(f,f) + f(v)(\mu - L[f]).$$

Such a $\mu$ might not exist as $L[f]$ is in general unbounded. Even if we assume $f$ to have compact support, choosing a $\mu$ to satisfy (78) would lead to excessive numerical viscosity. Therefore, one seeks a $\mu$ just big enough to keep $P$ positive. Various heuristics (see [43, 16]) exist to achieve that.

We begin by introducing a change of variables:

$$\tau := 1 - e^{-\mu t/k_n}, \qquad F(\tau, v) := f(t, v) e^{\mu t/k_n}. \tag{79}$$

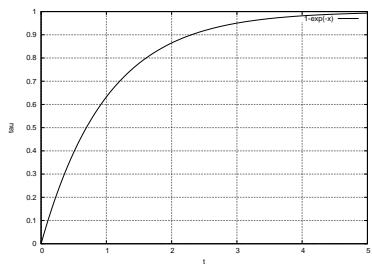

**Figure 9.** $\tau$ vs. $t$ in (79).

It is not hard to see that (79) changes (77) into

$$\frac{\partial F}{\partial \tau} = \frac{1}{\mu} P(F, F). \tag{80}$$

Formally expanding $F$ in a power series in $t$ yields

$$F(\tau, v) = \sum_{k=0}^{\infty} \tau^k \varphi_k(v), \tag{81}$$

and we immediately see that $\varphi_0(v) = f(x, v, 0)$. Substituting (81) into (80) gives a recurrence formula for $\varphi_{k+1}$:

$$\varphi_{k+1}(v) = \frac{1}{k+1} \sum_{l=0}^{k} \frac{1}{\mu} P(\varphi_l, \varphi_{k-l}). \tag{82}$$

It turns out that $\varphi_k \to f_\infty$ as $k \to \infty$, where $f_\infty$ is the Maxwellian. Truncating (81) and adding the "limit term" of the series back on gives

$$f^{n+1}(v) = (1-\tau) \sum_{k=0}^{m} \tau^k \varphi_k^n(v) + \tau^{m+1} f_\infty(v), \tag{83}$$

where $\varphi_k^n$ denotes $\varphi_k$ as in (82), but is calculated from $f^n$, the value of $f$ at time step $n$. Note that the factor $(1-\tau)$ is added to ensure that

$$1 = (1-\tau) \sum_{k=0}^{m} \tau^k + \tau^{m+1}.$$

Other choices of the coefficients of $\varphi_k$ and $f_\infty$ are conceivable, see [16]. We examine the properties of this scheme:

**Proposition 9.** *(cobbled together from [16] and [43])*

*The scheme (83) has the following properties:*

1. *It is consistent with (77) of order m.*

2. *For $m \in \{1, 2\}$, these schemes are A-stable.*

3. *If $f^n$ is nonnegative, then $f^{n+1}$ is also nonnegative.*

4. *Certain integrals are preserved: If $g$ is a function such that*

$$\int_{\mathbb{R}^3} P(f,f) g \, dv = \mu \int_{\mathbb{R}^3} f g \, dv,$$

*then*

$$\int_{\mathbb{R}^3} f^{n+1} g \, dv = \int_{\mathbb{R}^3} f^n g \, dv.$$

5. *The scheme preserves the asymptotic behavior of (77):*

$$\lim_{\mu \Delta t / k_n \to \infty} f^{n+1}(v) = f_\infty(v).$$

# 5 Summary

In this report we have surveyed many of the deterministic numerical algorithms used to solve the collisional Boltzmann equation (1) that have been proposed and implemented in the literature. Recall that the main difficulties for solving the Boltzmann equation arise from the following:

- Dimensionality of the unknown: the density $f$ is a function of 6+1 space+time variables.

- Developing consistent time evolution methods to deal with fundamentally disparate spatial algorithms for the transport and the collision terms.

- Computation of the five-fold integral collision term from equations (2), (3), and (4).

- Maintaining positivity of the density $f$ and accurately resolving high gradients.

Individually, none of the tasks above is trivial, yet solving the Boltzmann equation requires addressing all of them simultaneously.

In Section 2, we discussed various methods for fusing fundamentally different (spatial) discretizations with high-order time-stepping methods. Probably the most common, and the easiest to implement, are the splitting methods. The ubiquitous Strang splitting [48] is second-order accurate. An additional second-order accurate splitting method tailored for the Boltzmann equation is presented by Ohwada [38], and higher-order splitting methods can be found in [13].

An alternative to splitting methods are high-order implicit-explicit (IMEX) methods, which directly address coupling an explicit and an implicit solver. These methods are introduced by Ascher [2] and are elaborated upon by Pareschi [40]. They have the advantage of using a direct error analysis to obtain high-order temporal accuracy. They even have Strong Stability-Preserving properties under certain conditions.

In Section 3 we presented different methods for discretizing the transport form of the Boltzmann equation, or the Vlasov equation (26). This equation is essentially a hyperbolic conservation law, and there is a large body of literature which discusses solving such equations. Focusing on balancing the Boltzmann equation's need for low computational cost due to dimensionality, reasonable resolution of shock waves, and positivity of the solution, the methods we present can be classified as either fully discrete or semi-discrete specifically designed for the Vlasov equation. We first presented the fully-discrete flux balance method [14], which follows the characteristics and requires a reconstruction process to advance the scheme. This reconstruction process can be linear, cubic [18], or can employ (W)ENO methods [10]. A second choice is the semi-Lagrangian method. This method follows characteristics like the flux-balance method, and thus also requires a reconstruction. This reconstruction can be done with Lagrange or Hermite interpolations, and (W)ENO. Good surveys for the semi-Lagrangian methods can be found in [47] and [17].

We concluded our discussion of the transport term with a sample of some miscellaneous methods including finite-difference [1], spectral [32], and discontinuous Galerkin [20] discretizations.

Finally, Section 4 was devoted to surveying the deterministic methods for discretizing the collision term. A linearization can be obtained via the Bhatnagar-Gross-Krook approximation [49], which allows for much less painful evaluation of the collision term. However, if one wishes to solve the full Boltzmann equation, then one is forced to deal with the five-fold integral, which naively requires $O(N^8)$ operations per timestep, per point in space. Many efforts have been made to make this computation more tractable. The general modus operandi is to first truncate the generally infinite-domain velocity space, and then to make use of the FFT to decrease computational cost. Of course, with truncation, aliasing is a concern– Proposition 4, which gives sufficient conditions to avoid aliasing, partially addresses this issue. The first major reduction in cost is made by Ibragimov [29], down to $O(N^6)$ for certain collision kernel models ($O(N^6 \log N)$ for general kernel forms) with a quadratic order of convergence. Pareschi [41] improves on this by keeping the $O(N^6)$ operation count for certain collision kernels and obtains spectral convergence. Finally, the computationally fastest method for computing the collision operator in this survey is given by Mouhot in [35]. The method discretizes solid scattering angles with $M^2$ degrees of freedom and computes the operator with spectral accuracy in $O(M^2 N^3 \log N)$ computational time, albeit for a somewhat restricted set of collision kernels. To finish off, we addressed the issue of timestepping for the collision operator in the last part of Section 4.

Deterministic methods for the Boltzmann equation offer compelling advantages for problems where high accuracy and low noise is sought. For problems with low dimensionality, deterministic methods have already become a viable alternative to the more common Monte Carlo approaches. Furthermore, their accuracy makes them valuable as a source of validation for larger-scale randomized codes. We foresee a bright future for deterministic computation in this area and expect that forthcoming advances, both in faster machines and better algorithms, will continue to work in their favor.

**Acknowledgments.** We would like to thank Prof. Chi-Wang Shu for his advice and many helpful discussions. We are much indebted to Ishani Roy, who took a great amount of time and wrote quite a bit of text to educate us on the intricacies of applying the Boltzmann to the simulation of semiconductors. We would further like to thank our classmates in AM281 for the inspiring environment in this seminar course, for which this report was originally written.